\documentclass[11pt,a4paper,oneside,english]{article}
\usepackage[T1]{fontenc}
\usepackage[utf8]{inputenc}
\usepackage[english]{babel}
\usepackage{latexsym}
\usepackage{midpage}
\usepackage{newlfont}
\usepackage{color}
\usepackage{array}
\usepackage{amsmath} 
\usepackage{amsfonts} 
\usepackage{amsthm} 
\usepackage{amssymb}
\usepackage{setspace}
\usepackage[hidelinks]{hyperref} 
\usepackage{url}
\usepackage[normalem]{ulem}
\usepackage{cleveref}
\usepackage{graphicx} 
\usepackage{caption}
\usepackage{subcaption} 
\usepackage{mathdots}
\usepackage{longtable}
\usepackage{booktabs}
\usepackage[margin=1in,left=1.04in,right=1.04in,includefoot]{geometry}
\usepackage{rotating}
\usepackage{listings}
\usepackage{arydshln}
\usepackage{tikz}
\usepackage{color}
\usepackage{bbm}
\usepackage{bm}
\usepackage{comment}

\definecolor{darkgreen}{rgb}{0, .5, 0}

\theoremstyle{plain} 
\newtheorem{thrm}{Theorem}[section] 
\newtheorem{cor}[thrm]{Corollary} 
\newtheorem{lem}[thrm]{Lemma} 
\newtheorem{remark}[thrm]{Remark} 
\newtheorem{prop}[thrm]{Proposition} 
\theoremstyle{definition} 
\newtheorem{defn}{Definition}[section] 
\theoremstyle{remark}

\numberwithin{equation}{section}

\providecommand{\abs}[1]{\ensuremath{\left\lvert#1\right\rvert}}
\newcommand{\E}{\mathbb{E}} 
\newcommand{\N}{\mathbb{N}} 
\newcommand{\p}{\mathbb{P}} 
\newcommand{\Q}{\mathbb{Q}} 
\newcommand{\R}{\mathbb{R}} 
\newcommand{\F}{\mathcal{F}} 
\newcommand{\semi}{\mathcal{S}} 

\newcommand{\W}{\mathbb{W}} 
\newcommand{\q}{\mathcal{Q}} 
\newcommand{\h}{\mathcal{H}} 

\newcommand{\tg}{\tilde{\gamma}}

\newcommand{\norm}[1]{\left\lVert#1\right\rVert}

\newcommand{\Log}{\mathrm{Log}}
\newcommand{\Exp}{\mathrm{Exp}}

\author{Nils Detering\thanks{\textsc{Heinrich Heine University Düsseldorf, Faculty of Mathematics and Natural Sciences, Mathematics Institute}; \textbf{nils.detering@hhu.de}} 
	\and Silvia Lavagnini\thanks{\textsc{Department of Data Science and Analytics, BI Norwegian Business School, 0484 Oslo, Norway}; \textbf{silvia.lavagnini@bi.no}}}

\title{
A class of locally state-dependent models for forward curves
}

\begin{document}
\maketitle

\begin{abstract}
We present a dynamic model for forward curves within the Heath-Jarrow-Morton framework under the Musiela parametrization. The forward curves take values in a function space $\h$, and their dynamics follows a stochastic partial differential equation with state-dependent coefficients. In particular, the coefficients are defined through point-wise operating maps on $\h$, resulting in a locally state-dependent structure. We first explore conditions under which these point-wise operators are well defined on $\h$. Next, we determine conditions to ensure that the resulting coefficient functions satisfy local growth and Lipschitz properties, so to guarantee the existence and uniqueness of mild solutions. The proposed model captures the behavior of the entire forward curve through a single equation, yet retains remarkable simplicity. Notably, we demonstrate that certain one-dimensional projections of the model are Markovian and satisfy a one-dimensional stochastic differential equation. This connects our Hilbert-space approach to well established models for forward contracts with fixed delivery times, for which existing formulas and numerical techniques can be applied. This link allows us to examine also conditions for maintaining positivity of the solutions. As concrete examples, we analyze Hilbert-space valued variants of an exponential model and of a constant elasticity of variance model.
\end{abstract}

\section{Introduction}
We introduce a class of Hilbert-space valued state-dependent models for forward curves in commodity and interest rate markets. A forward curve is usually represented by a function $T\mapsto F(t,T)$, where $F(t,T)$ denotes the forward price at time $t\ge 0$ for a contract with instantaneous delivery in $T\ge t$. In some commodity markets, however, the contracts do not deliver at a precise time $T$, but rather over a time interval of the form $[T_1, T_2]$, with $T_2>T_1$, hence the price $F(t,T)$ is not observed directly. This is the case, e.g., in electricity and natural gas markets. However, it is still convenient for these markets to model the (fictitious) forward curve $F(t, \cdot)$ with instantaneous delivery. Then, the price at time $t \leq T_1$ of the forward with delivery period $[T_1, T_2]$ is obtained as the weighted integral of $F(t, \cdot)$ over the interval $[T_1, T_2]$. 

We focus in this article on the instantaneous forward curve $\{F(t, \cdot)\}_{t\ge 0}$, covering in this manner both commodity and interest rate markets. In particular, this approach guarantees that suitable no-arbitrage conditions arising from the possible overlap of contracts (e.g., contracts delivering over months and contracts delivering over quarters) are easily satisfied; see \cite{fredbook} for more details. Under the Musiela parametrization \cite{musiela} with $g_t(x):=F(t,T)$ for $x:=T-t$, one can then ensure that the time $t$ instantaneous forward curve $g_t$ takes values in a Hilbert space $\h$ of functions from $[0,\infty )$ to $\mathbb{R}$. Specifically, this function space can be chosen so that the forward curves have some empirically observed properties, and so that the dynamic behavior of $\{ g_t \}_{t\ge 0}$ can be described by a stochastic partial differential equation (SPDE).

A rather natural choice for $\h$ is the so-called Filipovi\'c space, which is the space of (weakly) differentiable functions whose weighted squared (weak) derivative is integrable. This space was first introduced in \cite{filipovic} in the context of Heath-Jarrow-Morton (HJM) term-structure modeling of bonds, and many of its mathematical properties have been derived there. Among other desirable features, the space reflects, for example, the well-known property that forward curves tend to flatten for $x\rightarrow \infty$, since less information is available for delivery times far in the future. We also refer to the recent book \cite{fredkrubook} with additional properties of $\h$, including an analysis of its algebraic structure and results related to energy markets. See further \cite{Benth2021,BDG2021,BDK2022,fredkru1,fredkru2,benthsgarra,clew, cox22} and \cite{jaber24} for additional applications. We also mention that modifications of the space to twice (weakly) differentiable functions are suggested in \cite{filipovic22b, filipovic22a} and \cite{FILWIL18} in the context of term-structure models, while in \cite{fredkru3} the space is extended to complex-valued functions. 

We use $\h$ as the state space for our forward curve model. For a function $\psi:\R\to\R$, we then start our analysis by studying point-wise operating maps $\Psi:\h\to\h$, which we define for any $f\in\h$ by $\Psi(f)(x) := \psi(f(x))$, for $x\in\R_+^0$. These maps act locally on the curve and are therefore particularly simple. However, they are interesting from a theoretical perspective and still offer the flexibility to specify a state-dependent model. 
In a second step, we construct time-dependent specifications for $\Psi$. The class of functions obtained can be used to specify the drift and the diffusion coefficient of the SPDE modeling the dynamics of the forward curves. 
More specifically, we consider diffusion terms that are multiplicative operators with time- and state-dependent kernel. Namely, for every $t\in\R_+^0$ and $f\in\h$, we consider diffusion terms of the form $\mathcal{M}_{\Psi(t,f)}:\h\to\h$ with point-wise operating kernel $\Psi(t,f)$ so that $\mathcal{M}_{\Psi(t,f)}(h) = \Psi(t,f)\cdot h$, for any $h\in\h$. 

In order to obtain existence and uniqueness of (local) mild solutions from \cite{T2012}, we investigate sufficient conditions on $\psi$ to ensure that $\Psi(t,f)$ and $\mathcal{M}_{\Psi(t,f)}$ are locally Lipschitz, locally bounded and/or of linear growth. Here the function space $\h$, which is chosen for its natural economic properties, poses some mathematical challenges. This is mainly because of the nature of the norm which involves a weighted integral of the squared (weak) derivative. For a function $f\in\h$, the norm allows one to control the overall smoothness, but makes it difficult to bound $f'(x)$ based on $\norm{f}$ for fixed $x\geq 0$. In fact, if the coefficient functions are specified based on the point-wise operating maps, they usually satisfy only a local Lipschitz condition (see Remark \ref{remark:globallip}), which makes the existence results from \cite{T2012} particularly relevant for our case. In \cite{T2012}, local Lipschitz and linear growth conditions, respectively local Lipschitz and local boundedness conditions, are indeed shown to be sufficient for the existence and uniqueness of global mild solutions, respectively local mild solutions. We mention that existence and uniqueness results for mild solutions of jump-diffusion SPDEs with path-dependent coefficients can be found in \cite{filipovic2010b} and \cite{filipovic2010}, and for specifications with L\'evy noise in \cite{marinelli}.  Our approach differs from the ones listed, in the sense that we pose conditions directly on the function $\psi$ that defines the point-wise operating map $\Psi$.

After establishing existence and uniqueness, we show that certain one-dimensional projections of the mild solution of the SPDE are Markovian and solve a one-dimensional stochastic differential equation (SDE). This links the Hilbert-space approach to early statements of the HJM framework \cite{0543c6d5-020c-3e84-82d9-926d734be64f}, where the dynamics of $\{F(t,T)\}_{t\geq 0}$ is directly specified by a separate SDE for each maturity $T\ge t$, leading to an uncountable family of SDEs for the entire market (see also, for example, \cite[Sec. 23]{bjork2009arbitrage}). From this perspective, our results unify these two modeling approaches. 

For many applications, one is interested in (strictly-) positive curves. We therefore derive sufficient conditions on the coefficients of the SPDE to ensure that the resulting solutions remain within the convex cone of (strictly-) positive functions in $\h$. Notice that positivity preserving properties for mild solutions of SPDEs in the Filipovi\'c space have been considered in \cite{filipovic2010}, where both drift and diffusion coefficients are assumed to be globally Lipschitz and bounded. Extensions include \cite{T2022} and \cite{T2017}, where some of the technical assumptions are relaxed. Due to our point-wise specification, we are able to exploit the above mentioned link of the SPDE with one-dimensional SDEs, and to derive direct conditions for the positivity preserving property on the functions defining the coefficients. Moreover, we study restrictions on the kernel function in order to obtain invertibility of the multiplicative operator. This is a desirable property which allows one to change from the physical measure to an equivalent measure as, for example, the risk-neutral measure under which the discounted traded instruments are martingales. We demonstrate this in the follow-up paper \cite{DeteringLavagnini} for the case of commodity and energy markets, where we also perform a calibration to market data. 

In the last part of the paper, as two concrete examples, we consider an exponential class of models, extending the theory of geometric models in finite dimension to Hilbert-space valued models, and a Hilbert-space valued counterpart of a constant elasticity of variance (CEV) model. In particular, we show how the theoretical results in this paper can be applied to these two cases. 

The rest of the article is organized as follows: in Section \ref{section1} we introduce the setting in which we shall work and the SPDE for the forward curve. In Section \ref{sec:spde_coeff} we define a class of point-wise operating SPDE coefficients, ensure that they are well defined, and study growth and Lipschitz properties. In Section \ref{sec:prop} we analyze invertibility, positivity and the one-dimensional projection for our model. Finally, in Section \ref{sec:spec} we present the exponential and the CEV specifications, and in Section \ref{sec:conclusions} we conclude with some remarks and discuss possible extensions.

\section{The forward curve model}\label{section1}
Let $F(t, T)$ be the price at time $t\le T$ for a forward contract with maturity at time $T$. We define $F(t, T)= F(T, T)$ for any $t> T$. In some markets, such as, for example, the electricity market, this contract is a hypothetical contract that is not traded but serves as a modeling device; see, e.g. \cite{fredbook}. Following Musiela \cite{musiela}, we introduce the parametrization $g(t,x):= F(t, t+x)$ where $x:= T-t$ is the \emph{time to delivery}. In the same spirit as \cite{Benth2021} and \cite{filipovic2010}, we can then understand the entire forward curve $\left\{g(t,x)\right\}_{x\ge 0}$ for any time $t\ge0$ as an element in a suitable function space, and state a dynamical model for it. In the following, let $\R_{+}^0$, respectively $\R_+$, be the set of positive and strictly-positive real numbers.

\subsection{The space of forward curves}
Following Filipovi{\'{c}} \cite{filipovic}, we introduce a space $\h:= \h(\R_{+}^0)$ of functions defined on $\R_{+}^0$: for a given continuous and non-decreasing function $w : \R_{+}^0\to [1,\infty)$ with $w(0)=1$, this is the Hilbert space of all absolutely continuous functions $f: \R_{+}^0 \to \R$ for which
\begin{equation*}
	\int_{\R_+^0}f'(x)^2w(x)dx < \infty .
\end{equation*}
In particular, with the inner product and the norm defined respectively by
\begin{equation*}
	\langle f_1,f_2\rangle := f_1(0)f_2(0)+\int_{\R_+^0}f_1'(x)f'_2(x)w(x)dx, \qquad \mbox{for } f_1, f_2\in\h,
\end{equation*}
and $\norm{f_1}^2:= \langle f_1,f_1\rangle$, this space is a separable Hilbert space (see \cite[Prop. 3.4]{BenthHalpha2017}). We shall assume that 
\begin{equation}
	\label{eq:alphaint}
\underline{w} :=\left( \int_{\R_+^0}w^{-1}(x)dx\right)^{1/2}<\infty, 
\end{equation}
which ensures sufficiently fast decay of the derivative of the functions in $\h$. One can then show that for a function $f\in \h$ the limit $\lim_{x\rightarrow \infty} f(x)$ exists (see \cite[Lemma 3.2]{fredkru1}). A typical example for the weight function is $w(x)= e^{c x}$ with $c >0$. We shall denote by  $\mathcal{L}(\h)$ the space of linear and bounded operators from $\h$ to itself, equipped with the usual operator norm, here denoted by $\norm{\cdot}_{\mathrm{op}}$. Moreover, we use the notation $\langle\cdot, \cdot\rangle_\R$ for the scalar product in $\R$.

Following \cite[Sec. 4.1.1, Def. 4.2]{sdebook}, we now introduce an $\h$-valued Wiener process $\W :=\{\W_t\}_{t\ge 0}$ with a compact, self-adjoint and strictly-positive definite covariance operator $\q\in \mathcal{L}(\h)$. In particular, there exists an orthonormal basis $\{e_j\}_{j}$ of $\h$ and a bounded sequence $\{\ell_j\}_{j}$ of strictly-positive real numbers such that 
\begin{equation}
	\label{eq:Qexpansion}
	\q f = \sum_{j} \ell_j \langle f, e_j \rangle_{} e_j, \qquad \mbox{for any } f\in \h,
\end{equation}
namely, the ${\ell_j}$'s are the eigenvalues of $\q$, and each $e_j$ is an eigenvector corresponding to ${\ell_j}$. We assume that $\mathrm{tr}(\q) = \sum_j {\ell_j} < \infty$, namely $\q$ is a trace-class operator. 

We then introduce the space $\h^\q:= \q^{1/2}\left(\h\right)$, which, equipped with the inner product
\begin{equation*}
	\langle f_1, f_2\rangle_{\q} := \langle \q^{-1/2}f_1, \q^{-1/2}f_2\rangle,  \qquad \mbox{with } f_1,f_2\in \h^\q,
\end{equation*}
is a separable Hilbert space with orthonormal basis $\{\sqrt{{\ell_j}}e_j\}_{j}$. We further denote by $ \mathcal{L}_{\mathrm{HS}}\left( \h^\q, \h\right)$ the space of Hilbert-Schmidt operators from $\h^\q$ into $\h$, which, endowed with the Hilbert-Schmidt norm
\begin{equation*}
	\norm{\mathcal{A}}_{\mathrm{HS}} := \sqrt{\sum_j {\ell_j} \norm{\mathcal{A}e_j}^2}, \qquad \mbox{for } \mathcal{A}\in \mathcal{L}_{\mathrm{HS}}\left( \h^\q, \h\right),
\end{equation*}
is itself another separable Hilbert space.

While we model the entire forward curve, ultimately one may be interested in the value of the forward curve at certain time points (times to maturity) or within certain time intervals for the case of forwards with delivery period. For this reason, for any $x\in\R_{+}^0$ we introduce the point evaluation map $\delta_x: \h \to \R$ as $\delta_xf := f(x)$, for $f\in\h$. By \cite[Lemma 3.8]{BenthHalpha2017}, the map $\delta_x$ is a bounded linear operator. In particular, under assumption \eqref{eq:alphaint}, the positive constant $K_{\delta}$, defined via $K_{\delta}^2 := 2\max(1, \underline{w}^2)$, is a uniform bound for the operator norm of $\delta_\cdot$, namely
\begin{equation}\label{norm:bound}
\norm{\delta_x}_{\mathrm{op}} = \sup_{\norm{f}=1}|\delta_x f| \le  K_{\delta},
\end{equation}
for every $x\in\R_{+}^0$, where $\norm{\cdot}_{\mathrm{op}}$ denotes the operator norm in $\mathcal{L}(\h,\R)$. Moreover, from \cite[Lemma 3.1]{fredkru1} the dual operator of $\delta_x$ is 
\begin{equation*}
	\begin{aligned}
	\delta_x^*: \R &\to \h\\
	c&\mapsto \left\{y \mapsto c\left(1 + \int_0^{x\land y} w^{-1}(u)du\right)\right\},
	\end{aligned}
\end{equation*}
with norm bounded by
\begin{equation*}
	\norm{\delta_x^*(c)}^2 = 
	c^2\left(1+ \int_0^x w^{-1}(y)dy\right) \le c^2\left(1+\underline{w}^2\right)\le c^2\,2\max(1, \underline{w}^2) = c^2K_\delta^2,
\end{equation*}
for all $x\in\R_{+}^0$ and $c\in\R$.

\subsection{The forward curves dynamics}
We now specify a dynamics for the forward rates $\{g(t,\cdot) \}_{t\geq 0}$.  In the related literature, it is common to state the equation for the forward curve directly under a pricing measure $\Q$. This leads to an equation without drift in the case of commodity markets, see e.g. \cite{Benth2021},  and with a particular drift in the case of interest rate markets, see e.g. \cite{filipovic2010}. We state the dynamics of $g$ under the physical probability measure $\p$. This allows one to use statistical inference on the observed forward curve to set up a model that is empirically validated. In our forthcoming paper \cite{DeteringLavagnini}, we then study equivalent measures under which tradable one-dimensional projections are martingales in the context of commodity and energy markets. 

Let $\left(\Omega, \F, \{\F_t\}_{t\geq 0},\mathbb{P}\right)$ be a filtered probability space satisfying the usual conditions, and let a Wiener process $\{ \W_t \}_{t\geq 0}$ with covariance kernel $\q$ be defined on this probability space. We further assume that $\mathrm{tr}(\q)<\infty$. For $g_t := g(t, \cdot)$, we consider the following stochastic partial differential equation (SPDE)
\begin{equation}
	\label{spde}
	dg_t = \partial _x g_t dt +  \alpha(t, g_t) dt + \sigma (t,g_t) d\W_t,
\end{equation} 
where $\partial_x$ is the generator for the shift-semigroup $\{\mathcal{S}_t\}_{t\ge 0}$ given by $\mathcal{S}_tf(x)=f(t+x)$, for any $t,x\in \R_{+}^0$ and $f\in \h$. Moreover, the mapping
$\alpha : \R_{+}^0\times \h \rightarrow \h$ is a state-dependent drift term, and the operator $\sigma : \R_{+}^0\times \h\rightarrow \mathcal{L}_{\mathrm{HS}}(\h^\q,\h)$ is a state-dependent diffusion term. 

Because the $C_0$-semigroup $\{\mathcal{S}_t\}_{t\ge 0}$  is quasi-contractive by \cite[Prop. 4.4]{BDK2022}, one can show, under Lipschitz and growth (respectively local boundedness) conditions for $\alpha$ and $\sigma$, that a (local) mild solution to \eqref{spde} exists in terms of the semigroup. But before stating existence and uniqueness results, we recall some definitions that will be practical for our setting. In the following, let $\left(\mathcal{K}_1, \norm{\cdot}_{\mathcal{K}_1}\right)$ and $\left(\mathcal{K}_2, \norm{\cdot}_{\mathcal{K}_2}\right)$ be two Banach spaces. Let further $\mathcal{K} \subseteq \mathcal{K}_1$. 
\begin{defn}[Locally Lipschitz]
	\label{def1}
	We say that a map $b: \R_{+}^0\times \mathcal{K} \to \mathcal{K}_2$ is \emph{locally Lipschitz} if for each $n\in\N$ there exists a non-decreasing function $L_n:\R_{+}^0\to \R_{+}^0$ such that 
	\begin{equation*}
		\norm{b(t, f_1)-b(t, f_2)}_{\mathcal{K}_2} \le L_n(t)\norm{f_1-f_2}_{\mathcal{K}_1}
	\end{equation*}
	for all $t\in \R_{+}^0$ and for all $f_1, f_2\in \mathcal{K}$ such that $\norm{f_1}_{\mathcal{K}_1}, \norm{f_2}_{\mathcal{K}_1}\le n$. Moreover, if there exists a non-decreasing function $L:\R_{+}^0\to \R_{+}^0$ such that $L_n$ in the equation above can be replaced by $L$ for every $n\in\N$, then we simply say that $b$ is \emph{Lipschitz}.
\end{defn}

\begin{defn}[Linear growth]
	\label{def2}
	We say that a map $b: \R_{+}^0\times \mathcal{K} \to \mathcal{K}_2$ satisfies the \emph{linear growth} condition if there is a non-decreasing function $G:\R_{+}^0\to \R_{+}^0$ such that 
	\begin{equation*}
		\norm{b(t, f)}_{\mathcal{K}_2} \le G(t)\left(1+\norm{f}_{\mathcal{K}_1}\right)
	\end{equation*}
	for all $t\in \R_{+}^0$ and all $f\in \mathcal{K}$.
\end{defn}

\begin{defn}[Locally bounded]
	\label{def3}
	We say that a map $b: \R_{+}^0\times \mathcal{K} \to \mathcal{K}_2$ is \emph{locally bounded} if for each $n\in\N$ there is a non-decreasing function $B_n:\R_{+}^0\to \R_{+}^0$ such that 
	\begin{equation*}
		\norm{b(t, f)}_{\mathcal{K}_2} \le B_n(t)
	\end{equation*}
	for all $t\in \R_{+}^0$ and for all $f\in \mathcal{K}$ such that $\norm{f}_{\mathcal{K}_1}\le n$.
\end{defn}

It is easy to notice that linear growth in Definition \ref{def2} implies locally bounded in Definition \ref{def3}.  Moreover, locally Lipschitz together with the condition that $\sup_{s\in [0,t]} b(s,\bm{0}) <\infty $, where $\bm{0}$ is the zero vector, implies locally bounded. In our setting, $\bm{0}$ will always be the function constantly equal to zero. We also point out that for functions $b$ which do not depend on time, namely $b: \mathcal{K} \to \mathcal{K}_2$, the non-decreasing functions introduced in the three definitions above are in fact constant, namely $L_n(t)=L_n$, $G(t)=G$ and $B_n(t)=B_n$, for all $t\in\R_{+}^0$.

Next, we provide the definition of (local) mild solutions.

\begin{defn}[Mild solution]
	A \emph{local mild solution} to \eqref{spde} is a pair $(g, \tau)$, where $\tau$ is a stopping time such that $\tau > 0$ almost surely, and $g:\Omega \times [0, \tau)\to \h$ is a measurable adapted process with continuous trajectories such that
	\begin{equation*}
		\E\left[ \int_0^t  \mathcal{S}_{t-s}\sigma(s, g_s)\q \left( \mathcal{S}_{t-s}\sigma(s, g_s)\right)^* ds \right]< \infty,
	\end{equation*} and, for any initial condition $g_0= g(0, \cdot )\in \h$, the identity
	\begin{equation}
		\label{mildsol}
		g_{t} = \mathcal{S}_{t} g _0 + \int_0^t \mathcal{S}_{t-s}\alpha(s, g_s)ds +  \int_0^{t}\mathcal{S}_{t-s}\sigma (s, g_s) d\W_s
	\end{equation}
	holds almost surely for any stopping time $t$ such that $t< \tau$ almost surely. 
	Moreover, a local mild solution $(g, \tau)$ is a \emph{global mild solution}, or simply a \emph{mild solution}, if $\tau = +\infty$ almost surely.
\end{defn}

The existence and uniqueness of mild solutions are guaranteed by \cite[Thm. 4.5]{T2012}. In particular, it holds that:
	\begin{enumerate}
		\item If $\alpha$ and $\sigma$ are locally Lipschitz and locally bounded, then there exists a unique local mild solutions to \eqref{spde};
		\item If $\alpha$ and $\sigma$ are locally Lipschitz and of linear growth, then there exists a unique global mild solutions to \eqref{spde}.
	\end{enumerate}

In Section \ref{sec:spde_coeff} we shall study sufficient conditions on our specific $\alpha : \R_{+}^0\times \h \rightarrow \h$ and $\sigma : \R_{+}^0\times \h\rightarrow \mathcal{L}_{\mathrm{HS}}(\h^\q,\h)$ to guarantee the conditions above, and hence existence and uniqueness of (local) mild solutions. Regarding the diffusion term $\sigma$, however, one can show that, without loss of generality, we can focus on operators in $\mathcal{L}(\h)$ instead of $\mathcal{L}_{\mathrm{HS}}(\h^\q,\h)$. For $\sigma(t,f)\in  \mathcal{L}_{\mathrm{HS}}(\h^\q,\h)$ with $t\in \R_{+}^0$ and $f\in\h$, we have indeed
\begin{equation*}
    \norm{\sigma(t,f)}_{\mathrm{HS}}^2 =  \sum_j {\ell_j} \norm{\sigma(t,f)e_j}^2 \le \sum_j {\ell_j}  \norm{\sigma(t,f)}_{\mathrm{op}}^2 = \norm{\sigma(t,f)}_{\mathrm{op}}^2\sum_j {\ell_j}  =\norm{\sigma(t,f)}_{\mathrm{op}}^2\mathrm{tr}(\q),
\end{equation*}
where in our setting $\mathrm{tr}(\q)<\infty$. Hence, since $\h^\q \subset \h$, it follows from $\sigma(t,f) \in \mathcal{L} (\h)$ that $\sigma(t,f)\vert_{\h^\q} \in \mathcal{L}_{\mathrm{HS}}(\h^\q,\h)$. This further implies that if $\sigma : \R_{+}^0\times \h\rightarrow \mathcal{L} (\h)$ is (locally) Lipschitz and/or locally bounded and/or of linear growth, then so is $\sigma\vert_{\h^\q} : \R_{+}^0\times \h^\q\rightarrow \mathcal{L}_{\mathrm{HS}}(\h^\q,\h)$. Hence we shall focus on diffusion terms of the form $\sigma : \R_{+}^0\times \h\rightarrow \mathcal{L}(\h)$.

We point out that the result in \cite[Thm. 4.5]{T2012} about existence and uniqueness of global mild solutions under local Lipschitz and linear growth conditions is particularly relevant for our setting, as it turns out that a global Lipschitz property for point-wise operating maps from $\h$ to $\h$ rarely holds, as we shall see.

\section{The SPDE coefficients}
\label{sec:spde_coeff}
We study in this section a specification for the drift term $\alpha$ and for the diffusion term $\sigma$ in the SPDE \eqref{spde}, which depend both on $t$ and on the current forward curve $g_t$.
In particular, to specify $\alpha$, we need to consider functions that map from $\R_{+}^0\times\h$ into $\h$. While, in order to specify $\sigma$,  we need a sufficiently rich class of functions that map from $\R_{+}^0\times\h$ into $\mathcal{L} (\h)$. 
 
In particular, for specifying the diffusion term, we observe that every function $h \in\h$ can be used to introduce a multiplicative operator $\mathcal{M}_h$. More precisely, the \emph{multiplicative operator $\mathcal{M}_h$ with kernel} $h \in\h$ is defined by
\begin{align*}
	\mathcal{M}_h: \h &\to \h\\
	f &\mapsto  h f.
\end{align*} 
It has been shown in \cite[Thm. 4.17]{fredkru1} that, under assumption \eqref{eq:alphaint}, $\mathcal{M}_h$ is well defined for every $h\in \h$, namely $h f \in \h$ for all $f\in \h$, and it is obviously linear. Moreover, it can be shown to be bounded with 
\begin{equation}\label{norm:mult:operator}
	\norm{\mathcal{M}_h}_{\text{op}}\leq K_{\mathcal{M}} \norm{h},
\end{equation}
where $K_{\mathcal{M}} := \sqrt{5+4\underline{w}^2}$. Hence $\mathcal{M}_h \in \mathcal{L}(\h)$ for all $h\in \h$. It seems then reasonable to specify $\sigma$ as a multiplicative operator with time-dependent and state-dependent kernel function of the form $h(t, g_t)$, with $h: \R_+^0\times \h  \to \h$. 

We then start in Section \ref{sec:pointwisemaps} by investigating maps $\Psi: \h  \to \h$ which can be used to specify $\alpha$ and the kernel for a multiplicative operator $\mathcal{M}_{\Psi}$. In particular, for any $f\in\h$ and $x\in \R_+^0$, we restrict to maps $\Psi$ which are such that $\delta_x\Psi(f) = \Psi( f)(x)= \psi(f(x))$ for some function $\psi : \R \rightarrow \R$. We say in the following that the map $\Psi$ \emph{operates point-wise} on functions $f$ in the space $\h$. These maps are simple but yet flexible enough to allow the resulting model to be calibrated to market data. We therefore start by first investigating suitable conditions on $\psi $ which ensure that the map 
\begin{equation}
	\label{eq:Psi}
	\begin{aligned}
		\Psi:  \h & \to \h \\
		f&\mapsto  \{x\mapsto \psi(f(x))\}
	\end{aligned}
\end{equation}
is well defined. By studying this class of maps, we obtain state-dependent functions to be used both for the drift coefficient $\alpha$ and for the kernel function $\Psi$ of the multiplicative operator $\mathcal{M}_{\Psi}$. In a second step, in Section \ref{sec:conditions} we then consider functions $\Psi$ which also depend on time, so to specify coefficients of the form $\alpha: \R_{+}^0\times\h\to\h$ and $\sigma=\mathcal{M}_{\Psi}:\R_{+}^0\times\h\to\mathcal{L} (\h)$ which may account for seasonality and a maturity-dependent volatility structure.

For some applications, we are interested in restricting the forward model to (strictly-) positive curves. We therefore denote by $\h_+^0$ and $\h_+$, respectively, the convex cone of positive and strictly positive functions:
\begin{align*}
&\h_+^0 := \{f \in \h \, : \, f(x) \ge 0 \mbox{ for every } x \in \R_+^0 \},\\
&\h_+ := \{f \in \h \, : \, f(x) > 0 \mbox{ for every } x \in \R_+^0 \}.
\end{align*}
In particular, we shall derive in Section \ref{sec:positivity} sufficient conditions on the coefficients of the SPDE  \eqref{spde} to ensure solutions in $\h_+^0$ and in $\h_+$.  Moreover, when working in $\h_+^0$, we shall restrict the domain of a point-wise operator $\Psi$ in \eqref{eq:Psi} to $\h_+^0 \subset \h$, hence $\psi$ maps from $\R_+^0$ to $\R$. Similarly, when working in $\h_+$,  the domain of $\Psi$ is $ \h_+\subset \h$ and $\psi$ maps from  $\R_+$ to $\R$. To ease the notation in the following results, we introduce generic sets $\mathcal{X}$ with $\mathcal{X}\in \{\R, \R_{+}^0, \R_{+}\}$ and $\mathcal{Y}$ with $\mathcal{Y}\in \{\h, \h_{+}^0, \h_{+}\}$. Without further mentioning, $\mathcal{X}= \R$ always implies that $\mathcal{Y}= \h$, while $\mathcal{X}= \R^0_+$ implies that $\mathcal{Y}= \h_+^0$, and $\mathcal{X}= \R_+$ implies that $\mathcal{Y}= \h_+$.

\subsection{Point-wise operating maps}
\label{sec:pointwisemaps}
We start by first investigating conditions on $\psi$ as a map from $\mathcal{X} \rightarrow \R$, with $\mathcal{X}\in \{\R, \R_{+}^0, \R_{+}\}$, ensuring that $\Psi (f) \in \mathcal{Y}$ for every $f\in \mathcal{Y}$, with $\mathcal{Y}\in \{\h, \h_{+}^0, \h_{+}\}$. The distinction between the different cases will be emphasized in all the results.  

\begin{prop}
\label{prop:mapping:into:h} 
Let $\psi: \mathcal{X} \rightarrow  \mathbb{R} $. The following statements hold:
\begin{enumerate}
\item $\mathcal{X} = \R$: If $\psi $ is differentiable with $\psi'$ continuous, then $\Psi(f)\in\h$ for every $f\in\h$.

\item $\mathcal{X} \in \{\R_{+}^0, \R_{+}\}$: If $\psi$ is differentiable in $\R_+$ with $\psi'$ continuous and such that  $\lim_{y\rightarrow 0} \psi'( y) <\infty $, then $\Psi(f)\in\h$ for every $f\in\mathcal{Y}$.
\end{enumerate}
\end{prop}
\begin{proof}
We start with \emph{1.}. Since $f\in\h$ is weakly differentiable and $\psi: \mathbb{R} \rightarrow  \mathbb{R} $ is continuously differentiable, the map $x\mapsto \psi(f(x))$ is also weakly differentiable by \cite[p. 215, Cor. 8.11]{Brezis2011}. In particular, one obtains that
\begin{align}
	\label{norm:pointwise:psi}
	\norm{\Psi( f) }^2&= \norm{\psi( f(\cdot))}^2= \psi(f(0))^2  + \int_0^{\infty}   \left( \psi'( f(x))\right)^2 f' (x)^2w (x)dx \notag \\
	&\leq  \psi(f(0))^2 + \int_0^{\infty}  \sup_{x\in\R^0_+} \left|\psi'( f(x))\right|^2  f' (x)^2w (x)dx.
\end{align} 
Since $|f(x)| \leq  K_{\delta} \norm{f}$, we notice that $f(x)\in [- K_{\delta} \norm{f},  K_{\delta} \norm{f}]$ for all $x\in\R_{+}^0$. By the continuity assumption on $y\mapsto \psi'( y)$ and because the interval $[- K_{\delta} \norm{f},  K_{\delta} \norm{f}]$ is compact, it follows  that
\begin{equation}\label{eq:C_f}
C_f :=\sup_{x\in \R_+^0} \abs{\psi'(f(x))} \leq   \sup_{y\in [- K_{\delta} \norm{f},  K_{\delta} \norm{f}]} \abs{\psi'(y)} <\infty.
\end{equation} 
Equation \eqref{eq:C_f} together with \eqref{norm:pointwise:psi} implies that $\norm{\Psi(f) }^2\leq \psi( f(0))^2 + \norm{f}^2C_{f}^2<\infty$. Thus $\Psi( f) \in \h$, which proves statement \emph{1.}.

For statement \emph{2.}, let us first assume $\mathcal{X}=\R_{+}$ and let $f\in\h_+$. By the assumption $\lim_{y\rightarrow 0} \psi'( y) <\infty $, it follows that $\sup_{y\in (0,K]} \abs{\psi'(y)} <\infty$ for any $K>0$, and in particular $\sup_{x\in \R_+} \abs{\psi'(f(x))} \leq \sup_{y\in (0,K_{\delta} \norm{f}]} \abs{\psi'(y)} <\infty$. With a similar line of reasoning as in \emph{1.}, we then conclude that $\Psi(f) \in \h$.
For the case $\mathcal{X}=\R_{+}^0$, we notice that $\lim_{y\rightarrow 0} \psi'( y) <\infty $ implies that $\psi$ is right differentiable in $0$, as it can be seen by the mean value theorem, and that  $\sup_{y\in [0,K]} \abs{\psi'(y)}<\infty$. Hence, as above, it follows that $\Psi(f) \in \h$ for $f\in\h_+^0$.
\end{proof}

In the last result, we required the function $\psi$ to be differentiable. However, the elements in $\h$ only need to be absolutely continuous. It is therefore natural to ask whether the function $\psi:\R\to \R$ may have (a few) points where the derivative is actually not defined. One might expect that in such a situation the result in Proposition \ref{prop:mapping:into:h} remains true as long as the set $\mathcal{P}_{\psi}:=\{ y | \psi ' (y) \text{ is undefined} \}$ has measure zero. However, this is actually not the case. Indeed, for $y\in \mathcal{P}_{\psi}$ and $f\in\h$, the set $f^{-1} (\{y \} )$ might have positive measure, which makes the analysis a bit more delicate. 

The following proposition shows that, under a topological conditions on $\mathcal{P}_{\psi}$, we can still obtain $\Psi(f) \in \h$ for every $f\in\h$. In the statement, the term {\em discrete set} refers to a set whose elements are isolated points. We recall that any discrete subset of $\R$ is at most countable.

\begin{prop}
\label{prop:non:diff}
Let $\psi: \mathcal{X} \rightarrow  \mathbb{R}$. The following statements hold:
\begin{enumerate}
\item $\mathcal{X}=\mathbb{R}$:
 Let the function $\psi$ be Lipschitz continuous. If $\psi$ is differentiable with bounded derivative except for a discrete set $\mathcal{P}_{\psi}$ of points, then $\Psi(f) \in \h$. 
 \item $\mathcal{X} \in \{\R_{+}^0, \R_{+}\}$: Let the function $\psi$ be Lipschitz continuous. If $\psi$ is differentiable in $\mathbb{R}_+$ with bounded derivative except for a discrete set $\mathcal{P}_{\psi}$ of points, then $\Psi(f) \in \h$.
\end{enumerate}
\end{prop}
\begin{proof}
We show the result for $\mathcal{X}=\R$, but the other cases can be treated similarly. Because $\psi $ is  Lipschitz continuous, it follows from \cite{leoni} that $\Psi(f) $ is an absolutely continuous function. In particular, the following chain rule applies:
	\begin{equation*}
		\left(\psi \circ f\right)' = \left(\psi'\circ f\right) f',
	\end{equation*}
	where, for $x\ge0$, the product $(\psi'\circ f)(x) f'(x)$ must be interpreted to be zero whenever $f'(x)=0$, independently of whether $(\psi'\circ f)(x)$ is defined. Let $U_f:= \{x\in\R_{+}: \, f'(x)\ne0\}$. We then have that
	\begin{align*}
		\norm{\psi\circ f}^2 &= \left(\left(\psi \circ f\right)(0)\right)^2 + \int_{\R_{+}}	\left(\left(\psi \circ f\right)'(x)\right)^2w(x)dx\notag\\
		& = \left(\left(\psi \circ f\right)(0)\right)^2 + \int_{U_f}	\left((\psi'\circ f)(x) f'(x)\right)^2w(x)dx. 
	\end{align*}
	The first term in the equation above is clearly bounded. To bound the second term, we first define the set $\mathcal{D}_f:=\{x \vert f(x) \in \mathcal{P}_{\psi} \}$ related to the possible kinks of $\psi$.
		If we can show that $\mu ( U_f  \cap \mathcal{D}_f ) =0$, where $\mu$ is the Lebesgue measure, then we obtain that 
		\begin{equation*}
			\int_{U_f}	\left((\psi'\circ f)(x) f'(x)\right)^2w(x)dx =  \int_{U_f\setminus \mathcal{D}_f }	\left((\psi'\circ f)(x) f'(x)\right)^2w(x)dx.
		\end{equation*}
		On the set $U_f\setminus \mathcal{D}_f$ we then have that $\psi'\circ f$ is well defined and, in particular, bounded, which allows us to conclude that 
		\begin{align*}
			\int_{U_f}	\left((\psi'\circ f)(x) f'(x)\right)^2w(x)dx 
			&\leq C  \int_{U_f \setminus \mathcal{D}_f}	\left(f'(x)\right)^2w(x)dx <\infty,
		\end{align*}
		from which it follows that $\Psi(t, f) \in \h$. 
		
		To show that, in fact, $\mu ( \mathcal{U}_f  \cap \mathcal{D}_f ) =0$, we need the following topological argument. For each $z \in  \mathcal{U}_f  \setminus \mathcal{D}_f$, we know that $f'(z)$ exists and that $f'(z)\neq 0$. Because $\mathcal{P}_{\psi}$ is discrete, for $\varepsilon>0$ we can select an open interval $\mathcal{O}_z:= (z-\varepsilon , z + \varepsilon)$ which is such that $\mathcal{D}_f \cap \mathcal{O}_z = \{z \}$. If follows by the choice of $\mathcal{O}_z$ that  $\mathcal{U}_f  \cap \mathcal{D}_f \subset \cup_{z \in\mathcal{U}_f  \cap \mathcal{D}_f } \mathcal{O}_z$. Because $\mathbb{R}$ is a hereditarily Lindel\"of space (see for example \cite[Sec. 3.8]{engelking1989general}), it follows that there exists a countable subcover, i.e. a subset $Z \subset \mathcal{U}_f  \cap \mathcal{D}_f$, such that $\mathcal{U}_f  \cap \mathcal{D}_f \subset \cup_{z \in Z} \mathcal{O}_z$. By construction, we have $\mu (\mathcal{O}_z \cap \mathcal{D}_f)=0$. Since $Z$ is countable, it then follows that $\mu ( \mathcal{U}_f  \cap \mathcal{D}_f ) \leq \mu ( \cup_{z \in Z} \mathcal{O}_z) =0$.
\end{proof}
\begin{remark}
Both Proposition \ref{prop:mapping:into:h} and Proposition \ref{prop:non:diff} threat the case $\psi:\mathcal{X} \to \R$ and $\Psi (f) \in \h$ with $\mathcal{X}\in \{\R, \R_{+}^0, \R_{+}\}$. However, if $\psi$ is either $\psi:\mathcal{X} \to \R_+^0$ or $\psi:\mathcal{X} \to \R_+$, then the same statements trivially imply $\Psi (f) \in \h_+^0$, respectively $\Psi (f) \in \h_+$.
\end{remark}

\subsection{Growth and Lipschitz conditions}
\label{sec:conditions}
We now consider time-dependent specifications for $\Psi$, which will serve to define the coefficients of the SPDE \eqref{spde}, namely we consider $\Psi$  of the form
\begin{equation}
	\label{eq:Psi:t}
	\begin{aligned}
		\Psi: \R^0_{+}\times \h & \to \h \\
		(t, f)&\mapsto  \{x\mapsto \psi(t, f(x))\},
	\end{aligned}
\end{equation}
for $\psi: \R_{+}^0\times\mathcal{X} \rightarrow \mathbb{R}$,  with $\mathcal{X}\in \{\R, \R_{+}^0, \R_{+}\}$. 
We start with the following result which provides sufficient conditions on $\psi$ to ensure that $\Psi(t, f)$ is locally bounded. 
\begin{prop}
\label{prop:locallybounded} 
Let $\psi: \R_{+}^0\times\mathcal{X} \rightarrow  \mathbb{R} $. The following statements hold:
\begin{enumerate}
\item $\mathcal{X} = \R$: Let $\psi (t, \cdot ) $ be differentiable for all $t\in\R_{+}^0$. If further $\psi$ and $\partial_y\psi $ are continuous in $\R^0_+ \times \R$, then $(t,f) \mapsto \Psi(t, f)$ is locally bounded.
\item $\mathcal{X} = \R_{+}^0$: Let $\psi (t, \cdot )$ be differentiable in $\R_+$ and such that  $\lim_{y\rightarrow 0}  \partial_y\psi(t, y) <\infty $ for all $t\in\R_{+}^0$. Let further $\psi$ be continuous in $\R_+^0 \times \R_+^0$, $\partial_y\psi $ be continuous in $\R_+^0 \times \R_+$, and the function $b_1 : \R_{+}^0\rightarrow \R$ defined by $b_1(t):=\lim_{y\rightarrow 0}  \partial_y\psi(t, y)$ be continuous. Then $(t,f) \mapsto \Psi(t, f)$ is locally bounded.

\item $\mathcal{X} = \R_{+}$: Let $\psi (t, \cdot )$ be differentiable in $\R_+$ and such that  $\lim_{y\rightarrow 0}  \partial_y\psi(t, y) <\infty $ for all $t\in\R_{+}^0$. Let further the functions $\psi$ and $\partial_y\psi $ be continuous in $\R_+^0 \times \R_+$, and the functions $b_1 : \R_{+}^0\rightarrow \R$ and $b_2 : \R_{+}^0\rightarrow \R$, defined respectively by $b_1(t):=\lim_{y\rightarrow 0}  \partial_y\psi(t, y)$ and $b_2(t):=\lim_{y\rightarrow 0}  \psi(t, y)$, both be continuous. Then $(t,f) \mapsto \Psi(t, f)$ is locally bounded.
\end{enumerate}
\end{prop}

\begin{proof}
We start with \emph{1.}. Since $\psi$ fulfills the assumptions of case \emph{1.} of Proposition \ref{prop:mapping:into:h} for all $t\in\R_{+}^0$, we know that $\Psi(t, f)\in \h$ for all $f\in\h$ and $t\in\R_{+}^0$. We then need to show that $(t,f) \mapsto \Psi(t, f)$ is locally bounded. With similar steps as in the proof of Proposition \ref{prop:mapping:into:h}, we obtain that 
\begin{align}
	\norm{\Psi(t, f) }^2 &\leq  \psi(t, f(0))^2 + \int_0^{\infty}  \sup_{x\in\R^0_+} \left| \partial_y \psi(t, f(x))\right|^2  f' (x)^2w (x)dx,\label{normfgamma3}
\end{align} 
for all $f\in\h$ and $t\in\R_{+}^0$. 
Since $|f(x)| \leq  K_{\delta} \norm{f}$, we observe again that $f(x)\in [- K_{\delta} \norm{f},  K_{\delta} \norm{f}]$ for all $x\in\R_{+}^0$. By the continuity assumption on $y\mapsto \partial_y\psi(t, y)$, it follows  that for every $t\in\R^0_+$
\begin{equation*}
C_n(t) :=\sup_{x\in \R_+^0, \norm{f} \leq n} \abs{\partial_y\psi(t, f(x))} \leq   \sup_{y\in A_n} \abs{\partial_y\psi(t, y)} <\infty ,
\end{equation*} 
with $A_n:=[-K_\delta n, K_\delta n]$.  We now show that the continuity of $(t,y)\mapsto \partial_y\psi(t, y)$ implies that $C_{n}(t)$ depends continuously on $t$. For this, we need to show that for a fixed $t\geq 0$, for any $\varepsilon>0$ there exists a neighborhood $T_t$ of $t$ such that 
\begin{equation*}
\abs{C_{n}(t)-C_{n}(s)} = \abs{\sup_{y\in A_n}\abs{\partial_y\psi(t, y)} -\sup_{y\in A_n}\abs{\partial_y\psi(s, y)} } \leq \varepsilon
\end{equation*}
for any $s\in T_t$. However, we notice that the inequality above is implied by
\begin{equation}\label{eps:ineq}
	\sup_{y\in A_n}\abs{\partial_y\psi(t, y) -\partial_y\psi(s, y) } \leq \varepsilon,
\end{equation}
hence we shall proceed to show \eqref{eps:ineq}. Let $t\geq 0$ and $\varepsilon >0$ be fixed. By the continuity of $(t,y)\mapsto \partial_y\psi(t, y)$, for every $y\in A_n$ there exists an open neighborhood $O_y$ of $(t,y)$ in $\R_+^0 \times A_n$ such that for every $(s,\bar{y}) \in O_y$ it holds that $\abs{\partial_y\psi(t, y) -\partial_y\psi(s, \bar{y}) } \leq \varepsilon$. Let the  neighborhood $O_y$ be of the form $O_y:=T_y\times F_y$, where $T_y$ is an open interval in $\R_+^0$ and $F_y$ is an open interval in $\R$. Clearly $\{t \} \times A_n \subset \cup_{y\in  A_n} O_y$. Then, by compactness, one can choose a finite sub-cover such that $\{t \} \times A_n \subset \cup_{i=1}^L O_{y_i}$ for some $L\in \mathbb{N}$. Let now $T_t:=\cap_{i=1}^L T_{y_i} $.  By the choice of $O_y$, equation \eqref{eps:ineq} holds for any $s \in T_t$, and $C_{n}(t)$ is therefore continuous in $t$. In particular, the continuity of $C_{n}(t)$ implies that $\widetilde{C}_{n}(t)  := \sup_{s\in[0,t]} C_{n}(s)<\infty $, and clearly $\widetilde{C}_{n}$ is a non-decreasing function in $t$. 

A similar argument can be used for the function $\psi$ to show that $D_{n}(t):= \sup_{y\in A_n} \abs{\psi(t, y)}$ depends continuously on $t$. We can then define $\widetilde{D}_{n}(t)  := \sup_{s\in[0,t]} D_{n}(s)<\infty $, which is also non-decreasing in $t$.

Combining everything with equation \eqref{normfgamma3}, we obtain that
\begin{equation}\label{eq:Mn}
	\norm{\Psi(t, f) }^2 \leq \widetilde{D}_n(t)^2 + \widetilde{C}_n(t)^2 \norm{f}^2\leq \widetilde{D}_n(t)^2 + \widetilde{C}_n(t)^2 n =:M_n(t)^2<\infty, \quad \mbox{for }\norm{f}\leq n.
\end{equation}
Clearly $M_n(t)$ is non decreasing in $t$ by construction. This implies that $(t,f) \mapsto \Psi(t, f)$ is locally bounded and completes the proof of part \emph{1.}. 

For part \emph{2.}, since $\psi (t, \cdot )$ fulfills the assumptions of case \emph{2.} of Proposition \ref{prop:mapping:into:h} for all $t\in\R_{+}^0$, it follows that $\Psi(t, f)\in \h$ for all $t\in\R_{+}^0$ and $f\in\h_+^0$. The continuity assumptions on $\psi$, $\partial_y\psi $ and $b_1$ allow one to extend the definition of $\partial_y\psi $ to $\R_{+}^0\times \R_{+}^0$ in a continuous way. The arguments from the case \emph{1.} above can now be repeated with the compact set $A_n:=[0, K_\delta n]$.

For part \emph{3.}, since again $\psi (t, \cdot )$ fulfills the assumptions of case \emph{2.} of Proposition \ref{prop:mapping:into:h} for all $t\in\R_{+}^0$, it follows that $\Psi(t, f)\in \h$ for all $t\in\R_{+}^0$ and $f\in\h_+$.   Moreover, the continuity assumptions on $\psi$, $\partial_y\psi $, $b_1$ and $b_2$ allow one to extend the definition of $\psi$ and $\partial_y\psi $ to $\R_{+}^0\times \R_{+}^0$ in a continuous way. Again, the argument from \emph{1.} can be repeated with the compact set $A_n:=[0, K_\delta n]$. This concludes the proof.

\end{proof}

\begin{remark}
We stress that the assumptions on the functions $b_1$ and $b_2$ in Proposition \ref{prop:locallybounded} are satisfied if $\psi$ and $\partial_y\psi $ are uniformly continuous in $\R_+ \times \R_+$. However, we will consider in Section \ref{sec:cev} examples where only the weaker assumptions posed directly on $b_1$ and $b_2$ hold. 
\end{remark}

We now derive sufficient conditions on $\psi$ which ensure linear growth for the map $\Psi$. Obviously, these conditions are stronger than the ones stated in Proposition \ref{prop:locallybounded}.

\begin{prop}
	\label{prop:lin:growth} 
Let the map $\psi: \R_{+}^0\times\mathcal{X} \rightarrow  \mathbb{R}$ be of linear growth and let $\sup_{(t,y) \in  \R_{+}^0\times\mathcal{X} } \partial_y \psi (t,y) < \infty$. Then $(t,f) \mapsto \Psi(t, f)$ is of linear growth as a map from $\R_{+}^0\times\mathcal{Y} $ to $\mathcal{H} $.
\end{prop}
\begin{proof}
Let $C:=\sup_{(t,y) \in  \R_{+}^0\times\mathcal{X} }\left| \partial_y \psi (t,y)\right| < \infty$ and let $G(t)$ be such that $\abs{\psi (t,y)} \leq G(t) (1 + \abs{y})$ for $(t,y) \in  \R_{+}^0\times\mathcal{X}$. We obtain that
\begin{align*}
	\norm{\Psi(t, f) }^2 &\leq  \psi(t, f(0))^2 + \int_0^{\infty}  \sup_{x\in\R^0_+} \left| \partial_y \psi(t, f(x))\right|^2  f' (x)^2w (x)dx \notag \\
	& \leq G(t)^2 (1+ f(0))^2  + C^2 \norm{f}^2  \notag\\
		& \leq G(t)^2 (1+ \norm{f})^2  + C^2 \norm{f}^2.
\end{align*} 
Combining the inequality above with $\sqrt{x+y}\leq \sqrt{x} +\sqrt{y}$ for $x,y\geq 0$, it follows that $\norm{\Psi(t, f) } \leq  ( G(t) +C ) (1+ \norm{f})$, from which we get that $(t,f) \mapsto \Psi(t, f)$ is of linear growth in all the three cases, namely for all $\mathcal{X}\in \{\R, \R_{+}^0, \R_{+}\}$.
\end{proof}
 \begin{remark}We further stress that the boundedness assumption on $ \partial_y \psi $ required in Proposition \ref{prop:lin:growth} is often fulfilled when $\psi$ satisfies the linear growth assumption. However, the simple (even globally bounded) example $\psi(x)=\sin(\exp(x))$ shows that this is not always the case, and the assumption is therefore needed. 
 \end{remark}

We now show that, under appropriate conditions on $\psi$, the map $\Psi$ is locally Lipschitz.

\begin{prop}
\label{prop:lipchitz}
Let $\psi: \R_{+}^0\times\mathcal{X} \rightarrow  \mathbb{R}$. The following statements hold:
\begin{enumerate}
\item $\mathcal{X} = \R$: Let $\psi $ be differentiable in the second variable and let $\partial_y \psi: \mathbb{R}_+^0 \times \mathbb{R} \rightarrow  \mathbb{R}$ be locally Lipschitz. Then $(t,f) \mapsto \Psi(t, f)$ is well defined and locally Lipschitz as a map from $\mathbb{R}_+^0 \times \h $ to $ \h$.

\item $\mathcal{X} \in \{\R_{+}^0, \R_{+}\}$: Let $\psi $ be differentiable in the second variable and let $\partial_y \psi: \mathbb{R}_+^0 \times \mathbb{R}_+ \rightarrow  \mathbb{R}$ be locally Lipschitz. Then $(t,f) \mapsto \Psi(t, f)$ is well defined and  locally Lipschitz as a map from $\mathbb{R}_+^0 \times\mathcal{Y} $ to $ \h$.
\end{enumerate}
\end{prop}
\begin{proof} In order to show that the map $(t,f) \mapsto \Psi(t, f)$  is well defined in the two cases, it is sufficient to  show that conditions \emph{1.} and \emph{2.} imply respectively conditions \emph{1.} and \emph{2.} of Proposition \ref{prop:mapping:into:h}. For \emph{1.}, this is obvious. For \emph{2.}, we first show that $\lim_{y\rightarrow 0}  \partial_y\psi(t, y) <\infty $ exists for every $t\in \R_+^0$. For this, let $\{y_n \}_{n\in \mathbb{N}}$ be a sequence in $\R_+$ with $\lim_{n\rightarrow \infty}y_n=0$. Then $\{y_n \}_{n\in \mathbb{N}}$ is bounded, and because the function $\partial_y\psi$ is locally Lipschitz, it follows that $\abs{\partial_y\psi(t, y_m)-\partial_y\psi(t, y_n)}\leq L(t) \abs{y_m - y_n}$ for some $L(t)$ and for all $m$ and $n$. This implies that $\partial_y\psi(t, y_n)$ is a Cauchy sequence and its limit exists. It follows that in all cases the map $(t,f) \mapsto \Psi(t, f)$ is well defined. 

For locally Lipschitz, let now $f_1,f_2\in \h$ with $\norm{f_1},\norm{f_2}\leq n$. By assumption, we know that $\psi$ has locally Lipschitz derivative in the second argument. It follows for all cases that $\psi$ itself is locally Lipschitz in its second argument. For case \emph{1.}, this follows directly from the mean value theorem. For case \emph{2.}, it follows from the mean value theorem and the existence of $\lim_{y\rightarrow 0} \psi(t, y)$, which can be shown in the same way as we showed that $\lim_{y\rightarrow 0}  \partial_y\psi(t, y) <\infty $ exists. In particular,  there exist non-decreasing functions $D_{n}:\R_{+}^0\to \R_{+}^0$ such that  
	\begin{equation*}
		| \psi (t, y_1)  - \psi (t, y_2) |\leq  D_n (t) | y_1 - y_2| ,
	\end{equation*}
	for $\abs{y_1}, \abs{y_2} \leq \left \lceil{n K_\delta}\right \rceil $. It follows that
	\begin{equation}
 \label{firststep}
	 | \psi (t, f_1(x))  - \psi (t, f_2 (x)) |\leq  D_n (t) | f_1(x) - f_2(x)| \le K_\delta D_n (t) \norm{f_1 - f_2}
	  	\end{equation}
	  for $\norm{f_1} \norm{f_2} \leq n$. 
	Let us now compute
	\begin{align}
		& \| \Psi (t,  f_1) - \Psi (t, f_2  ) \|^2\nonumber  \\
		&= (\psi (t, f_1 (0)) - \psi (t, f_2 (0))  )^2  + \int_0^{\infty} (f_1' (x)  \partial_y\psi (t, f_1(x) )   - f_2' (x) \partial_y\psi (t, f_2(x)) )^2 w (x) dx.\label{lip:int}
	\end{align}
	For the first term on the right-hand side, one obtains from equation \eqref{firststep} that
	$$(\psi (t, f_1 (0)) - \psi (t, f_2 (0)) )^2 \leq K_\delta^2 D_n (t)^2 \norm{f_1-f_2}^2 $$ for $\norm{f_1} \norm{f_2} \leq n$.
	For the second term we have
	\begin{align*}
		& \int_0^{\infty} (f_1' (x)  \partial_y\psi (t, f_1(x) )   - f_2' (x) \partial_y\psi (t, f_2(x)) )^2 w (x) dx \\
		&\leq  2 \int_0^{\infty} \left( f_1' (x)  - f_2' (x) \right)^2  \partial_y\psi (t, f_1(x))^2 w (x) dx \\
		& \quad +2\int_0^{\infty} f_2' (x)^2 (\partial_y\psi (t, f_1(x) ) -\partial_y\psi (t, f_2(x) ) )^2 w (x) dx .
	\end{align*}
	With similar steps as in the proof of Proposition~\ref{prop:locallybounded}, it follows that there exists $\tilde C_n$ such that $| \partial_y{\psi} (t, f_1(x) )| \leq \tilde C_n (t)$ for $\norm{f_1} \leq n$, which implies that the first term is bounded by $2  \tilde C_n(t)^2\| f_1 - f_2\|^2$ for $\norm{f_1}\leq n$. Moreover, similarly as for $\psi$ in \eqref{firststep}, there exists $\tilde{D}_n:\R_{+}^0\to \R_{+}^0$ such that 
	$| \partial_y{\psi} (t, f_1(x)) - \partial_y{\psi}(t, f_2(x)) | \leq K_\delta \tilde{D}_n (t) \norm{f_1-f_2}$. This allows us to bound the second term by 
	$$2K_\delta^2 \tilde{D}_n (t)^2 \norm{f_1-f_2}^2 \norm{f_2}^2 \leq 2K_\delta^2 \tilde{D}_n(t) ^2n^2 \norm{f_1-f_2}^2 , \quad \mbox{for }\norm{f_2}\leq n.$$
	Altogether, we obtain that 
	$$\| \Psi (t, f_1) - \Psi (t, f_2 ) \|^2_{} \leq  L_n(t)^2 \norm{f_1-f_2}^2 ,$$
	with Lipschitz constant ${L}_n(t):=(K_\delta^2 D_n (t)^2+2\tilde C_n (t)^2 + 2K_\delta^2 \tilde{D}_n (t)^2 n^2)^{1/2}$ for all $t\in\R_{+}^0$.
\end{proof}

\begin{remark}\label{remark:globallip}
From line \eqref{lip:int} in the proof of Proposition \ref{prop:lipchitz}, we notice that the map $\Psi$ is usually not Lipschitz even if $\partial_y\psi$ is Lipschitz. Assume for instance that there exists $L>0$ such that 
\begin{equation}\label{lip:cond:possible}
    \abs{f_1' (x)  \partial_y\psi (f_1(x) )   - f_2' (x) \partial_y\psi (f_2(x))} \leq L \abs{f'_1(x)-f'_2(x)} ,
\end{equation}
for any $f_1,f_2\in \h$ and all $x\geq 0$, where we dropped the dependency on $t$ for simplicity. In this case, the Lipschitz continuity of the map $f\mapsto \Psi (f)$ would certainly hold by \eqref{lip:int}. However, unless the function $\partial_y\psi$ is constant, there exist $y_1$ and $y_2$ such that $\partial_y\psi (y_1)\neq \partial_y\psi (y_2)$. Fix now any $\bar{x}\geq 0$ and choose $f_1$ and $f_2$ such that $f_1(\bar{x}) =y_1 $, $f_2(\bar{x}) =y_2 $ and $f'_1(\bar{x})=f'_2(\bar{x})\neq 0$, which is certainly possible in $\h$. It then follows that 
\begin{equation*}
\abs{f_1' (\bar{x})  \partial_y\psi (f_1(\bar{x}) )   - f_2' (\bar{x}) \partial_y\psi (f_2(\bar{x}))}
=\abs{f_1' (\bar{x}) \left(  \partial_y\psi (f_1(\bar{x}) )   - \partial_y\psi (f_2(\bar{x}))\right)}>0.
\end{equation*}
On the other hand, it follows from \eqref{lip:cond:possible} that 
\begin{equation*}
\abs{f_1' (\bar{x})  \partial_y\psi (f_1(\bar{x}) )   - f_2' (\bar{x}) \partial_y\psi (f_2(\bar{x}))} \leq L \abs{f'_1(\bar{x})-f'_2(\bar{x})}=0,
\end{equation*}
which is a contradiction. 
\end{remark}
An additional continuity assumption for $\psi$ on top of the one required in Proposition \ref{prop:lipchitz} allows us to obtain an alternative result to Proposition \ref{prop:locallybounded}  for the map $\Psi (t, f)$ to be locally bounded. In the following proposition we state conditions that ensure both the locally bounded and the locally Lipschitz properties.   

\begin{prop}
\label{prop:bounded}
Let $\psi: \R_{+}^0\times\mathcal{X} \rightarrow  \mathbb{R} $. The following statements hold:
\begin{enumerate}
\item $\mathcal{X} = \R$: Let $\psi $ be differentiable in the second variable and let the maps $\psi(t,y)$  and $\partial_y\psi(t,y)$ be continuous in $\R^0_+ \times \R$. If, in addition, $\partial_y \psi: \mathbb{R}_+^0 \times \mathbb{R} \rightarrow  \mathbb{R}$ is locally Lipschitz, then $(t,f) \mapsto \Psi(t, f)$ is locally bounded as a map from $\mathbb{R}_+^0 \times \h \rightarrow \h$.
\item $\mathcal{X} \in \{\R_{+}^0, \R_{+}\}$: Let $\psi $ be differentiable in the second variable and let the maps $\psi(t,y)$  and $\partial_y\psi(t,y)$ be continuous in $\R^0_+ \times \R_+$. If, in addition, $\partial_y \psi: \mathbb{R}_+^0 \times \mathbb{R}_+ \rightarrow  \mathbb{R}$ is locally Lipschitz, then $(t,f) \mapsto \Psi(t, f)$ is locally bounded as a map from $\mathbb{R}_+^0 \times\mathcal{Y} \rightarrow \h$.
\end{enumerate}
\end{prop}
\begin{proof}
We first show that the conditions of Proposition \ref{prop:locallybounded} are fulfilled. For \emph{1.}, this holds trivially, hence we only check \emph{2.}. Let first $\mathcal{X} = \R_{+}^0$, and let $b_1(t):=\lim_{y\rightarrow 0}  \partial_y\psi(t, y)$ for $t\in [0,\infty)$ which exists by the same reasoning used in Proposition \ref{prop:lipchitz}. We need to show that $b_1$ is continuous. For $t \in \R_+^0$ and $\varepsilon >0$,  we then show that $\abs{b_1(t)-b_1(s)}\leq \varepsilon$, for $s$ in a neighbourhood of $t$. Without loss of generality, we consider $s\leq t$. Let again $L(t)$ be the Lipschitz constant for the locally Lipschitz property of $\partial_y\psi$ for $y\leq 1$, and let $\bar{y}=\varepsilon/(3 L(t))\leq 1$, where we can reduce $\varepsilon$ if necessary. Now choose $\delta $ such that $\abs{\partial_y\psi(t, \bar{y}) -\partial_y\psi(s, \bar{y}) } \leq \varepsilon /3$ for $\abs{t-s}\leq \delta $. We further note that $\abs{\partial_y\psi(s,y) -\partial_y\psi(s, \bar{y}) } \leq \abs{\bar{y} - y} L(t) \leq \varepsilon /3$ for $y\in (0,\bar{y})$ and $s\leq t$. Therefore $\abs{ \partial_y\psi(s, \bar{y}) -  \lim_{y\rightarrow 0}  \partial_y\psi(s, y)}\leq \varepsilon /3 $. If follows that 
\begin{align*}
&\abs{b_1(s)-b_1(t)} = \abs{ \lim_{y\rightarrow 0}  \partial_y\psi(s, y) - \lim_{y\rightarrow 0}  \partial_y\psi(t, y)}\\\nonumber 
&\leq \abs{ \lim_{y\rightarrow 0}  \partial_y\psi(s, y) - \partial_y\psi(s, \bar{y})}+\abs{ \partial_y\psi(s, \bar{y}) - \partial_y\psi(t, \bar{y})}+\abs{ \partial_y\psi(t, \bar{y}) - \lim_{y\rightarrow 0}  \partial_y\psi(t, y)}\leq  \varepsilon .
\end{align*}
This shows that the function $b_1$ is continuous as required. For $\mathcal{X} = \R_{+}$, we observe that because $\partial_y \psi: \mathbb{R}_+^0 \times \mathbb{R}_+ \rightarrow  \mathbb{R}$ is locally Lipschitz, then also $ \psi: \mathbb{R}_+^0 \times \mathbb{R}_+ \rightarrow  \mathbb{R}$ is locally Lipschitz. This allows us to apply the same proof as for $\mathcal{X} = \R_{+}^0$ to show that the function $b_2: \R_+^0: \rightarrow \R$ defined by $b_2(t):=\lim_{y\rightarrow 0}  \psi(t, y)$ is well defined and continuous. The proof for the function $b_1: \R_+^0: \rightarrow \R$ defined by $b_1(t):=\lim_{y\rightarrow 0}  \partial_y\psi(t, y)$ is the same as in the case $\mathcal{X} = \R_{+}^0$. 

Hence we know by Proposition \ref{prop:locallybounded} that in all cases the map $(t,f) \mapsto \Psi(t, f)$ is locally bounded. 
\end{proof}

We end this section by considering a specification where the map $\Psi$ factorises into a time-dependent and a state-dependent component. Assume that there exist $\beta:\R_{+}^0\to \R$ and $\Phi:\mathcal{Y} \to \h$ with $\mathcal{Y} \in \{ \h, \h_+^0, \h_+ \} $ such that
\begin{equation}
	\label{eq:betaphi}
\Psi(t, f)= \beta(t)\Phi(f), \quad \mbox{for all } f\in \mathcal{Y} \mbox{ and } t\in\R_{+}^0.
\end{equation} 
The following results on $(t,f)\mapsto \Psi(t,f)$ when $\Phi (f)$ is of the form $\Phi (f)(x) =\phi (f(x)) $ for some real-valued function $\phi$ are then a consequence of Proposition \ref{prop:lin:growth}, Proposition \ref{prop:lipchitz} and Proposition \ref{prop:bounded}.
\begin{cor}\label{lem:specialcase:mult:structure}
	Let $\beta:\R_{+}^0\to \R$ be bounded,  
 and let $\phi: \mathcal{X} \rightarrow  \mathbb{R}$ be such that $\Phi (f)(x):=\phi (f(x))$. The following statements hold:
	\begin{enumerate}
		\item $\mathcal{X} = \R$: If $\phi $ is differentiable in $\R$ and $\phi'$ is locally Lipschitz continuous, then $\Psi (t, f) =\beta(t) \Phi (f) \in \h$ for every $f \in \h$. Moreover, the map $(t,f)\mapsto  \Psi (t, f)$ is locally bounded and locally Lipschitz. 
		
		\item  $\mathcal{X} \in\{\R_+^0, \R_+\}$: If $\phi $ is differentiable in $\R_+$ and $\phi' $ is locally Lipschitz continuous, then $\Psi (t, f) =\beta(t) \Phi (f) \in \h$ for every $f \in \mathcal{Y}$. Moreover, the map $(t,f)\mapsto \Psi (t, f)$ is locally bounded and locally Lipschitz.
	\end{enumerate}
\end{cor}
\begin{proof}
	For $\psi(t, f(\cdot)) = \beta(t)\phi(f(\cdot))$, these are special cases of Proposition \ref{prop:lipchitz} and Proposition \ref{prop:bounded}. 
\end{proof}

\begin{remark}
    We remark that, if $\beta(t) \in \R_+$ and $\phi(y) \in \R_+$, respectively $\beta(t) \in \R_+^0$ and $\phi(y) \in \R_+^0$, then $\Psi (t, f) \in \h_+$, respectively $\Psi (t, f) \in \h_+^0$.
\end{remark}

\begin{cor}\label{lem:specialcase:mult:structure:lingrowth}
	Let $\beta:\R_{+}^0\to \R$ be bounded, and let $\phi: \mathcal{X} \rightarrow  \mathbb{R}$ be such that $\Phi (f)(x):=\phi (f(x))$. If $\phi$ is of linear growth and $\sup_{y\in \mathcal{X}}\phi'(y)<\infty$, then the map $(t,f)\mapsto  \Psi (t, f)$ is of linear growth as a map from $\mathcal{Y}$ to $\h$.
\end{cor}
\begin{proof}
	For $\psi(t, f(\cdot)) = \beta(t)\phi(f(\cdot))$, this is a special cases of Proposition \ref{prop:lin:growth}. 
\end{proof}
\begin{remark}
\label{remark:psix}
It is also possible to consider $\beta:\R_{+}^0\to \mathcal{H}$. In this case, we would have $\delta_x \Psi (t, f) = \delta_x \beta(t) \Phi(f) = \beta(t)(x) \phi(f(x))$ for all $t,x\geq 0$ and $f\in\h$. While $\Phi(f)$ is just a special case of the maps studied before, the $\Psi(t,f)$ defined by $\delta_x \Psi (t, f) = \delta_x \beta(t) \Phi(f) = \beta(t)(x) \phi(f(x))$ is not. This is due to the extra dependence on $x$ which is not through $f(x)$. However, we point out that under the assumptions of Corollary \ref{lem:specialcase:mult:structure} on $\phi$, the map $\Psi (t, f) = \beta(t) \Phi(f)$ is well defined in $\h$ due to the Banach algebra structure of $\h$ (see \cite[Thm. 4.17]{fredkru1}). To obtain statements along the lines of Corollary \ref{lem:specialcase:mult:structure} and Corollary \ref{lem:specialcase:mult:structure:lingrowth} for $\Psi$ we would need to ask that $\sup_{s\in [ 0,t ]} \norm{\beta (s)} < \infty$ for all $t\geq 0$, which holds, in particular, if $\beta:\R_{+}^0\to \mathcal{H}$ is continuous. Finally, we point out that if $\beta(t) \in \h_+$ and $\phi(y) \in \R_+$, respectively $\beta(t) \in \h_+^0$ and $\phi(y) \in \R_+^0$, then $\Psi (t, f) \in \h_+$, respectively $\Psi (t, f) \in \h_+^0$.
\end{remark} 

 \subsection{The multiplicative diffusion operator}
 \label{sec:moltop}
Thanks to Proposition~\ref{prop:mapping:into:h}, Proposition \ref{prop:lin:growth}, Proposition \ref{prop:lipchitz}, and Proposition \ref{prop:bounded}  we can now state  sufficient conditions on $\psi$ to guarantee that the map $(t,f)\mapsto \mathcal{M}_{\Psi(t, f)}$ is locally bounded, of linear growth and locally Lipschitz.  This ensures, together with appropriate assumptions on the drift coefficient $\alpha$, the existence and uniqueness of (local) mild solutions to the SPDE \eqref{spde} with diffusion operator of the form $\sigma(t, f)= \mathcal{M}_{\Psi (t, f)}$, for $t\in\R_+^0$ and $f\in\h$.
\begin{lem}\label{lipchitz:op:norm}
	Let $\psi: \R_{+}^0\times\mathcal{X} \rightarrow  \mathbb{R}$. The following statements hold:
	
	\begin{enumerate}
		\item $\mathcal{X} = \R$: If $\psi $ fulfills Assumptions 1. of Proposition \ref{prop:lipchitz} and Assumptions 1. of Proposition \ref{prop:bounded}, then the map $ \mathcal{M}_{\Psi}:\mathbb{R}_+^0\times \h \to \mathcal{L}(\h)$ defined by $(t,f)\mapsto  \mathcal{M}_{\Psi (t, f)}$ is locally bounded and locally Lipschitz.
		\item $\mathcal{X} \in\{ \R_{+}^0, \R_{+}\}$: If $\psi $ fulfills Assumptions 2. of Proposition \ref{prop:lipchitz} and Assumptions 2. of Proposition \ref{prop:bounded}, then the map $\mathcal{M}_{\Psi}:\mathbb{R}_+^0\times \mathcal{Y} \to \mathcal{L}(\h)$ defined by $(t,f)\mapsto  \mathcal{M}_{\Psi (t, f)}$ is locally bounded and locally Lipschitz.
  \end{enumerate}
\end{lem}
\begin{proof}
	In all cases we know that the map $(t,f) \mapsto \Psi(t, f)$ is locally Lipschitz from Proposition \ref{prop:lipchitz} and  locally bounded from Proposition \ref{prop:bounded}. Let now $M_n (t)$ be the non-decreasing function defined in \eqref{eq:Mn} such that $\norm{\Psi(t, f)} \leq M_n (t)$ for $\norm{f}\leq n$ and $t\in\R_{+}^0$. Since the map $h\mapsto  \mathcal{M}_{h}$ is bounded, it follows that 
	$$ \norm{\mathcal{M}_{\Psi (t, f)}}_{\mathrm{op}} \leq K_{\mathcal{M}} \norm{\Psi (t, f) }  \leq  K_{\mathcal{M}}{M}_n(t),$$
	for $f\in\h$ with $\norm{f}\leq n$, which implies that the map $(t,f)\mapsto  \mathcal{M}_{\Psi (t, f)}$ is locally bounded. Similarly, since the map $h\mapsto  \mathcal{M}_{h}$ is Lipschitz continuous with Lipschitz constant $K_\mathcal{M}$, it follows that 
	$$ \norm{\mathcal{M}_{\Psi (t, f_1)} - \mathcal{M}_{\Psi (t, f_2)}}_{\mathrm{op}} \leq K_{\mathcal{M}} \norm{\Psi (t, f_1) - \Psi (t, f_2)}  \leq  K_{\mathcal{M}}{L}_n(t) \norm{f_1-f_2},$$
	for $f_1,f_2$ with $\norm{f_1},\norm{f_2}\leq n$, where ${L}_n(t)$ is the local Lipschitz bound for $(t,f) \mapsto \Psi(t, f)$ found in the proof of Proposition~\ref{prop:lipchitz}. This implies that the map $(t,f)\mapsto  \mathcal{M}_{\Psi (t, f)}$ is locally Lipschitz. 
\end{proof}

\begin{lem}\label{lingrowth:op:norm}
	Let $\psi: \R_{+}^0\times\mathcal{X} \rightarrow  \mathbb{R}$ satisfy the assumptions of Proposition \ref{prop:lin:growth}. Then the map $ \mathcal{M}_{\Psi}:\mathbb{R}_+^0\times \mathcal{Y} \to \mathcal{L}(\h)$ defined by $(t,f)\mapsto  \mathcal{M}_{\Psi (t, f)}$ is of linear growth.
\end{lem}
\begin{proof}
    From Proposition \ref{prop:lin:growth} we know that $(t,f)\mapsto  \Psi (t, f)$ is of linear growth, hence there exists a non-decreasing function $G:\R_+^0\to \R_+^0$ such that $\norm{\Psi(t, f) } \leq   G(t) (1+ \norm{f})$ for all $t\in\R_+^0$ and $f\in\mathcal{Y}$. Combining this with $\norm{\mathcal{M}_{\Psi (t, f)}}_{\mathrm{op}} \leq K_{\mathcal{M}} \norm{\Psi (t, f) }$, we obtain the result.
\end{proof}

We conclude the treatment on the diffusion operators with the following corollaries to Lemma \ref{lipchitz:op:norm} and Lemma \ref{lingrowth:op:norm} for the specification \eqref{eq:betaphi} where $\Psi$ is of the form $\Psi (t, f) = \beta(t)\Phi(f)$.
\begin{cor}\label{lem:specialcase:mult:structure2}
	Let $\beta:\R_{+}^0\to \R$ be bounded, and let $\phi: \mathcal{X} \rightarrow  \mathbb{R}$ be such that $\Phi (f)(x):=\phi (f(x))$. The following statements hold: 
	\begin{enumerate}
		\item $\mathcal{X} = \R$: If $\phi $ is differentiable in $\R$ and $\phi'$ is locally Lipschitz continuous, then $\mathcal{M}_{\Psi}:\mathbb{R}_+^0\times \h \to \mathcal{L}(\h)$ defined by $(t,f)\mapsto  \mathcal{M}_{\Psi (t, f)}$ is locally bounded and locally Lipschitz.
		
		\item $\mathcal{X} \in\{ \R_{+}^0, \R_{+}\}$: If $\phi $ is differentiable in $\R_+$ and $\phi' $ is locally Lipschitz continuous, then $\mathcal{M}_{\Psi}:\mathbb{R}_+^0\times \mathcal{Y} \to \mathcal{L}(\h)$ defined by $(t,f)\mapsto  \mathcal{M}_{\Psi (t, f)}$ is locally bounded and locally Lipschitz.
	\end{enumerate}
\end{cor}
\begin{proof}
	We first notice that $f\mapsto \Phi (f)$ is a special case of the map considered in Lemma \ref{lipchitz:op:norm}, where the function $\Psi$ is constant in time. It is then easy to see that the conditions posed on $\phi$ in \emph{1.} and \emph{2.} imply the conditions  \emph{1.} and \emph{2.} in Lemma \ref{lipchitz:op:norm}. This implies in all cases that $\Phi (f) \in \h$, and that the map $(t,f)\mapsto  \mathcal{M}_{\Phi ( f)}$ is locally bounded and locally Lipschitz.
 
	It also holds that $\mathcal{M}_{\Psi(t, f)}= \beta(t)\mathcal{M}_{\Phi(f)}$.
	Together with $\beta$ bounded, it then follows that $(t,f)\mapsto  \mathcal{M}_{\Psi (t, f)}$ is locally bounded as well.
	For the Lipschitz continuity, we similarly observe that 
	\begin{equation*}
			\norm{\mathcal{M}_{\Psi(t, f_1)}-\mathcal{M}_{\Psi(t, f_2)}}_{\mathrm{op}}  = |\beta(t)|\norm{\mathcal{M}_{\Phi(f_1)}-\mathcal{M}_{\Phi(f_2)}}_{\mathrm{op}}. 
	\end{equation*}
	Hence, using again that $\beta$ is bounded, it follows that  $(t,f)\mapsto  \mathcal{M}_{\Psi (t, f)}$ is locally Lipschitz.
\end{proof}

\begin{cor}\label{lem:specialcase:mult:structure3}
	Let $\beta:\R_{+}^0\to \R$ be bounded, and let $\phi: \mathcal{X} \rightarrow  \mathbb{R}$ be such that $\Phi (f)(x):=\phi (f(x))$ satisfies the assumptions of Corollary \ref{lem:specialcase:mult:structure:lingrowth}. Then the map $ \mathcal{M}_{\Psi}:\mathbb{R}_+^0\times \mathcal{Y} \to \mathcal{L}(\h)$ defined by $(t,f)\mapsto  \mathcal{M}_{\Psi (t, f)}$ is of linear growth. 
\end{cor}
\begin{proof}
    We first notice that $f\mapsto \Phi (f)$ is a special case of the map considered in Lemma \ref{lingrowth:op:norm}, where the function $\Psi$ is constant in time. This implies that the map $(t,f)\mapsto  \mathcal{M}_{\Phi ( f)}$ is of linear growth. Combining this with $\beta$ bounded, it then follows that $(t,f)\mapsto  \mathcal{M}_{\Psi (t, f)}$ is also of linear growth.
\end{proof}

\begin{remark}
\label{rmk:phi}
	The proofs of Corollary \ref{lem:specialcase:mult:structure2} and Corollary \ref{lem:specialcase:mult:structure3} show that for kernel functions of the form of \eqref{eq:betaphi}, under suitable conditions on the function $\beta$, it is sufficient to focus on the state-dependent component, since, indeed, if $\mathcal{M}_{\Phi(f)}\in\mathcal{L}(\h)$, then also $\mathcal{M}_{\Psi(t, f)}\in\mathcal{L}(\h)$, and if $ \mathcal{M}_{\Phi(f)}$ is locally Lipschitz and/or locally bounded and/or of linear growth, then also $\mathcal{M}_{\Psi(t, f)}$ is locally Lipschitz and/or locally bounded and/or of linear growth.
\end{remark}

 \begin{remark}
We point out that it is common in the literature to consider diffusion operators which are actually Hilbert-Schmidt operators, see, e.g.,  \cite{BenthHalpha2017}, \cite{fredkru1}, \cite{filipovic2010} and \cite{filipovic2010b}. However, the multiplicative operator $\mathcal{M}_{\Psi(t, f)}$ does not belong to this class in general. Indeed, in order to have $\mathcal{M}_{\Psi(t, f)}\in \mathcal{L}_{\mathrm{HS}}(\h, \h)$ we would need, by definition, that 
 \begin{equation*}
 \norm{\mathcal{M}_{\Psi(t, f)}}_{\mathrm{HS}}^2=\sum_{j} \norm{\mathcal{M}_{\Psi(t, f)} e_j}^2 <\infty, \quad \mbox{ for all } f\in\h,
 \end{equation*}
 where $\{e_j \}_{j} $ is an orthonormal basis of $\h$, see \cite[Appendix C]{sdebook}. One easily observes that this does not hold if $\Psi(t, f)$ is a constant function. A constant kernel function is obtained, for example, when $\Psi(t, f)= \beta(t)$ with $\beta:\R_+^0\to\R_+^0$ a continuous function, namely when the volatility operator does not depend on the time to maturity. In this case $\sum_{j} \norm{\mathcal{M}_{\Psi(t, f)} e_j}^2 =\left|\beta(t)\right|^2\sum_{j} \norm{e_j}^2  = \infty$.  Notice however, that $\mathcal{M}_{\Psi(t, f)}\in \mathcal{L}_{\mathrm{HS}}(\h^\q, \h)$. Indeed, in this case, $\norm{\mathcal{M}_{\Psi(t, f)}}_{\mathrm{HS}}^2=\sum_{j} \norm{\mathcal{M}_{\Psi(t, f)} \sqrt{\ell_j}e_j}^2= \left|\beta(t)\right|^2\sum_{j} \ell_j\norm{e_j}^2 = \left|\beta(t)\right|^2\sum_{j} \ell_j<\infty$.
 \end{remark}

\subsection{Invertibility of multiplicative operators}
Multiplicative operators are natural choices for the diffusion term in the SPDE \eqref{spde}. We now show that often, in fact, they are invertible, which allows, for example, to perform changes of measure. We shall use this result in Section \ref{sec:positivity} in order to prove positivity of the SPDE solutions.  

In particular, for $h\in\h$, we define the \emph{inverse operator} of $\mathcal{M}_{h}$ as the operator $\mathcal{M}_{h}^{-1}\in\mathcal{L}(\h)$ such that the compositions $\mathcal{M}_{h}^{-1}\circ \mathcal{M}_{h}$ and $\mathcal{M}_{h}\circ \mathcal{M}_{h}^{-1}$ correspond to the identity operator. To define the inverse operator, however, we need to restrict the kernel $h$ to a subset of $\h$. Hence, we first define the following two sets:
	\begin{align*}
		\h_>&:=\{ h\in \h \,|\, h(x)>0 \; \forall x\in\R^0_{+}  \mbox{ with } \lim_{x\rightarrow \infty}  h(x)>0  \}, \mbox{ and}\\ 
		\h_<&:=\{ h\in \h \,|\, h(x)<0 \; \forall x\in\R^0_{+} \mbox{ with }  \lim_{x\rightarrow \infty}  h(x)<0  \}.
	\end{align*}
We state the following auxiliary lemma.
\begin{lem}\label{bounded:fcts:in:h} 
Let $h \in \h_> \cup \h_<$. Then $\inf_{x\in \mathbb{R}^0_+}\abs{h(x)} > 0$.
\end{lem}
\begin{proof}
We first notice that, as a straightforward consequence of equation \eqref{norm:bound}, for $h \in \h$ we have $\sup_{x\in \mathbb{R}^0_+}\abs{h(x)} < \infty$. Let now $h\in \h_>$. By \cite[Lemma 3.2]{fredkru1} and assumption \eqref{eq:alphaint}, the limit $C:=\lim_{x\rightarrow \infty} h(x)>0$ exists. Then there exists $\bar{x}\in \R_{+}$ such that $h(x)\in [C/2,2C]$ for any $x\geq \bar{x}$. However, the function $h$ is continuous and therefore attains a minimum on the interval $[0,\bar{x}]$. Let $\bar{C}:=\inf_{x\in [0,\bar{x}]}{h(x)}$. Then $h(x)\geq \min \{\bar{C}, C/2 \}>0$ for all $x\in\R_{+}^0$, hence $\inf_{x\in \mathbb{R}_+^0}h(x) > 0$. The same argument can be applied for $h\in \h_<$ to the function $-h$, concluding the proof.
\end{proof}

\begin{remark}
\label{remark:infh}
    Even though Lemma \ref{bounded:fcts:in:h} ensures existence of $\inf_{x\in \mathbb{R}^0_+}\abs{h(x)} > 0$, it is not possible to construct a lower bound for $\inf_{x\in \mathbb{R}^0_+}\abs{h(x)}$ in terms of $\norm{h}$, for $h \in \h_> \cup \h_<$. In fact, one can show that for any $\varepsilon > 0$, there exist functions $h_1,h_2 \in \h_> \cup \h_<$ such that $\inf_{x\in \mathbb{R}^0_+}\abs{h_1(x)} =\inf_{x\in \mathbb{R}^0_+}\abs{h_2(x)}=\varepsilon$, but $\norm{h_1}=\varepsilon$ and $\norm{h_2}>\varepsilon$. For this, one can choose $h_1$ to be the function constantly equal to $\varepsilon$, and $h_2$, for example, to be $h_2 (x):=\varepsilon +x {\bm 1}_{x\leq y}+y{\bm 1}_{x> y}$, for $y> 0$. Then $\norm{h_2}=\left(\varepsilon^2 + \int_0^y w(x) dx\right)^{1/2}$. Since $w: \mathbb{R}_+^0 \rightarrow [1,\infty )$ is increasing and $w(0)=1$, it follows that $\norm{h_2}> \varepsilon$  for any $y>0$.
\end{remark}

We can now show the invertibility of the multiplicative operator. For this, for every $h\in\h_>\cup \h_<$, we introduce the notation $1/h$ for the map defined by
\begin{equation*}
	\begin{aligned}
		\frac{1}{h}: \R^0_{+} & \to \R \\
		x&\mapsto  \frac{1}{h(x)}.
	\end{aligned}
\end{equation*}
Notice that, for $h\in \h_>\cup \h_<$, the map $x\mapsto 1/h(x)$ is well defined.

\begin{prop}\label{prop:invertibility}
Let $h\in \h_>$, respectively $h\in \h_<$. Then $1/h\in \h_>$, respectively $1/h\in\h_<$. Moreover, $\mathcal{M}_{1/h}$ is the inverse operator of $\mathcal{M}_{h}$ for every $h\in \h_>\cup \h_<$.
\end{prop}
\begin{proof}
We know that for $h\in \h_>\cup \h_<$ the map $x\mapsto 1/h(x)$ is well defined. Moreover, since the map $y\mapsto 1/y$ is continuously differentiable in $(-\infty, 0 )  \cup ( 0,\infty )$ and $h\in \h_>\cup \h_<$, it follows that $1/h$ is weakly differentiable with derivative $-h'(x)/h(x)^2$ by \cite[p. 215, Cor. 8.11]{Brezis2011}. By Lemma~\ref{bounded:fcts:in:h}, it follows that $C:=\sup_{x\in\R^0_+ } 1/h(x) < \infty$.
We then obtain
\begin{align*}	
	\norm{\frac{1}{h}}^2_{}& =\left(\frac{1}{h(0)}\right)^2 + \int_0^{\infty} \frac{h'(x)^2}{h(x)^4}  w (x)dx \leq C^2  + C^4 \int_0^{\infty}  h'(x)^2 w (x)dx \le  C^2  + C^4 \norm{h}^2_{} <  \infty, \nonumber
\end{align*}	
	which implies that $1/h\in \h$. Clearly, then $1/h\in \h_>$ for $h\in \h_>$, and $1/h\in \h_<$ for $h\in \h_<$. The fact that $\mathcal{M}_{1/h}$ is the left and right inverse of $\mathcal{M}_{h}$ follows directly from the definition. 
\end{proof}

\section{Properties of forward curves}\label{sec:prop}
In the previous sections we studied state-dependent and point-wise operating functions $\Psi$ of the form \eqref{eq:Psi:t}. For our setting, such functions $\Psi$ are intended for specifying, for $t\in\R_{+}^0$ and $f\in\h$, both the drift term $\alpha(t,f)$ and the kernel function of the multiplicative diffusion operator $\sigma(t, f) = \mathcal{M}_{\Psi(t, f)}$, as studied in Section \ref{sec:moltop}. Specifically, for the drift term, we introduce a map $a:\R_+^0\times\R\to \R$ such that
\begin{equation}
	\label{eq:alpha}
	\begin{aligned}
		\alpha: \R^0_{+}\times \h & \to \h \\
		(t, f)&\mapsto  \{x\mapsto a(t, f(x))\},
	\end{aligned}
\end{equation}
or, equivalently, $\delta_x\alpha(t, f) = \alpha(t, f)(x)= a(t, f(x))$, for every $f \in \h$ and $t,x\in \R_{+}^0$, for which all the results in Section \ref{sec:pointwisemaps} and Section \ref{sec:conditions} apply.

\subsection{Dynamics of a forward with fixed delivery}
Working with state-depended point-wise operating coefficient functions and multiplicative operators has some great advantages. While the SPDE \eqref{spde} models the evolution of the entire forward curve $g_t$ with time $t\ge0$, when evaluating the curve in a point $x\in \R_+^0$, we obtain in fact a stochastic differential equation (SDE) for the fixed-delivery forward price. Remember that $\delta_x g_t = F(t, t+x) = F(t, T)$, for $T=t+x$. Hence the dynamics of a particular forward contract with maturity $T$ is given by an SDE. This property has the advantage that for pricing an option on the fixed-delivery forward $F(t, T)$ either analytic formulas or one-dimensional simulation schemes can be used. We make this formal in the next proposition.

\begin{prop}
	\label{prop:sde}
	Let the diffusion coefficient $\sigma$ be of the form $\sigma(t, f) = \mathcal{M}_{\Psi(t, f)}$, where $\Psi$ acts point-wise as in \eqref{eq:Psi:t} for some $\psi:\R_{+}^0\times \R\to \R$, and let also the drift coefficient $\alpha$ act point-wise as in \eqref{eq:alpha} for some $a:\R_{+}^0\times \R\to \R$. For $g_0 \in\h$, let both $\sigma$ and $\alpha$ satisfy locally Lipschitz and locally bounded conditions, respectively locally Lipschitz and linear growth conditions, and let $(g, \tau)$ be the corresponding (local) mild solution. Then, for every $T\ge 0$ there exists a one-dimensional Wiener process $W$ such that the  forward price $F(\cdot, T)$ follows the dynamics
	 	\begin{equation}\label{eq:sde_d}
	 	dF(t,T) = a (t ,F(t, T)) dt+ c_t(T) \psi(t, F(t, T)) dW_t, \qquad \mbox{for every }t< \tau,
	 \end{equation}
	 with $c_t(T):= \left(\sum_{j} {\ell_j}  \bigl(e_j(T-t) \bigr)^2\right)^{1/2}$. Here $\{e_j\}_{j}$ is an orthonormal basis of $\h$ and $\{{\ell_j}\}_{j}$ are the corresponding eigenvalues for the covariance operator $\q$ of $\W$, see equation \eqref{eq:Qexpansion}.
\end{prop}
\begin{proof}
		By \cite[Thm. 4.5]{T2012}, for any $g_0 \in\h$ there exists a unique (local) mild solution $(g, \tau)$ to the SPDE \eqref{spde}, where $\tau$ is a stopping time almost surely positive and $g$ is of the form
	\begin{equation*}
		g_{t} = \mathcal{S}_{t} g _0 + \int_0^{t}\mathcal{S}_{t-s}\alpha (s, g_s) ds+ \int_0^{t}\mathcal{S}_{t-s}\sigma (s, g_s) d\W_s,
	\end{equation*}
	where $\sigma (s, g_s) = \mathcal{M}_{\Psi (s, g_s)}$ and $t<\tau$. In particular, $\tau = \infty$ if the solution is global.  
Since $\delta_x \mathcal{S}_t = \delta_{t+x}$ for all $t,x\in\R_{+}^0$, it follows that 
	\begin{equation}
		F(t,T)=\delta_{T-t} g_{t} =  g _0 (T)+ \int_0^{t}\delta_{T-s}\alpha (s ,g_s) ds+ \int_0^{t}\delta_{T-s}\sigma (s ,g_s) d\W_s.\label{eq:Fpos}
	\end{equation}
Let now $\Gamma_s := \delta_{T-s}\sigma(s, g_s): \h \to \R$ with dual operator $\Gamma_s^*: \R \to \h$. One can then write
\begin{equation*}
	\Gamma_s\bigl(d\W_s\bigr) = \langle 1, \Gamma_s\bigl(d\W_s\bigr)\rangle_\R = \langle  \Gamma_s^*\left(1\right), d\W_s\rangle.
\end{equation*}
In particular, by definition of the Wiener process $\W$, the right-hand side of the equation above has Gaussian distribution with variance $\langle \q \Gamma_s^*\left(1\right), \Gamma_s^*\left(1\right)\rangle ds = \langle \Gamma_s \q \Gamma_s^*\left(1\right), 1\rangle_\R ds$. Following \cite[Thm. 2.1]{fredkru1}, this allows us to introduce a one-dimensional Wiener process $W$ such that 
\begin{equation*}
	\Gamma_sd\W_s = \left( \Gamma_s \q \Gamma_s^*\left(1\right)\right)^{1/2} dW_s.
\end{equation*}
For every $0\le t \le T$, we then rewrite the stochastic integral on the right-hand side of equation \eqref{eq:Fpos} as
\begin{equation*}
		\int_0^{t}\delta_{T-s}\sigma (s ,g_s) d\W_s 
		=\int_{0}^{t}\bigl( \Gamma_s \q \Gamma_s^*\left(1\right)\bigr)^{1/2}dW_s.
\end{equation*}
This tells us that the forward price \eqref{eq:Fpos}, which is driven by the $\h$-valued Wiener process $\W$ with covariance operator $\q$, can be represented as driven by a one-dimensional Wiener process. In particular, equation \eqref{eq:Fpos} becomes	\begin{equation}
		F(t,T)= g _0 (T)+ \int_0^{t}a (s ,g_s(T-s)) ds+ \int_0^{t}\Sigma(s, g_s) dW_s,\label{eq:Fpos2}
	\end{equation}
	where $\Sigma(s, g_s)^2 := \delta_{T-s}\sigma(s, g_s)\q\sigma(s, g_s)^*\delta_{T-s}^*(1)$. 
	
	Let us focus on $\Sigma(s, g_s)^2$ for $0\le s\le T$. From equation \eqref{eq:Qexpansion}, we can represent the covariance operator $\q$ acting on $\sigma(s, g_s)^*\delta_{T-s}^*(1)\in \h$ by
	\begin{align}
		\notag\q \sigma(s, g_s)^*\delta_{T-s}^*(1) & = \sum_{j} {\ell_j} \langle  \sigma(s, g_s)^*\delta_{T-s}^*\left(1\right), e_j \rangle \, e_j
		\\&= \sum_{j} {\ell_j} \langle  \delta_{T-s}^*\left(1\right),\sigma(s, g_s)e_j \rangle \, e_j = \sum_{j} {\ell_j}  \left(\delta_{T-s}\sigma(s, g_s)e_j \right)  e_j, \label{eq:QL}
	\end{align}
	where the last equality is due to the very definition of the dual operator $\delta_x^*$, namely $\langle  \delta_x^*\left(1\right),f\rangle =\langle 1,\delta_xf\rangle_\R  =   \delta_{x}f$, for every $f\in\h$ and $x\in\R_{+}^0$. From \eqref{eq:QL} we now get that
	\begin{align}
		\notag\Sigma(s, g_s)^2 &= \delta_{T-s}\sigma(s, g_s)\biggl(\sum_{j} {\ell_j}  \left(\delta_{T-s}\sigma(s, g_s)e_j \right)  e_j \biggr)\\
		& = \sum_{j} {\ell_j} \delta_{T-s}\sigma(s, g_s) \left(\delta_{T-s}\sigma(s, g_s)e_j \right)  e_j = \sum_{j} {\ell_j}  \bigl(\delta_{T-s}\sigma(s, g_s)e_j \bigr)^2, \label{eq:LQL}
	\end{align}
	where the steps above are justified by the linearity of $\delta_{x}:\h\to \R$ and $\sigma(s, g_s):\h\to\h$. From \eqref{eq:LQL} we then deduce that
	\begin{equation}
    \label{eq:LQL2}
        \begin{aligned}
		\Sigma(s, g_s)^2  &= \sum_{j} {\ell_j}  \bigl(\sigma(s, g_s)(T-s)\,e_j(T-s) \bigr)^2 = c_s(T)^2 \bigl(\sigma(s, g_s)(T-s) \bigr)^2 \\
		& =  c_s(T)^2\bigl(\Psi(s, g_s)(T-s) \bigr)^2 =  c_s(T)^2\bigl(\psi(s, g_s(T-s)) \bigr)^2 = c_s(T)^2 \bigl(\psi(s, F(s, T)) \bigr)^2, 
	\end{aligned}
    \end{equation}
	for $c_s(T):= \left(\sum_{j} {\ell_j}  \bigl(e_j(T-s) \bigr)^2\right)^{1/2}$. In particular, for any $T>0$, we have that $c_s(T)>0$  for all $0\le s\le T$, since all the ${\ell_j}$'s are strictly positive, and
	\begin{equation*}
		 c_s(T)^2=\sum_{j} {\ell_j}  \bigl(e_j(T-s) \bigr)^2 \le \sum_{j} {\ell_j}  \norm{e_j}^2\norm{\delta_x}^2_{\mathrm{op}} \le K^2_{\delta} \sum_{j} {\ell_j} < \infty,
	\end{equation*}
	since $\sum_{j}{\ell_j} < \infty$ and $\{e_j\}_j$ is an orthonormal basis. Hence $0 < c_s <\infty$  for all $0\le s\le T$.  For a fixed $T>0$, we can then rewrite equation \eqref{eq:Fpos2} by using \eqref{eq:LQL2} as follows:
	\begin{equation}\label{eq:sde}
		F(t,T) = F(0,T)+ \int_0^{t}a (s ,F(s, T)) ds+ \int_0^{t}c_s(T) \psi(s, F(s, T)) dW_s,
	\end{equation}
	which is a one-dimensional SDE corresponding to \eqref{eq:sde_d}, for all $t<\tau$.
\end{proof}

We emphasize that the Wiener process $W$ introduced in equation \eqref{eq:sde_d} for the dynamics of $F(\cdot, T)$ depends on the maturity time $T$, namely, for every $T\ge 0$ one has a different Wiener process for each fixed-delivery forward with delivery $T$. In the following proposition, we compute the covariance structure between the different forward contracts.

\begin{prop}
	Consider the same setting as in Proposition \ref{prop:sde}. Then, for every $T_1, T_2\ge 0$, there exist two one-dimensional Wiener processes $W^1$ and $W^2$ such that the forward price $F(\cdot, T_i)$, for $i=1,2$, follows the dynamics \eqref{eq:sde_d},
	and the correlation $\rho$ between $W^1$ and $W^2$ is
	\begin{equation*}
		  \rho =\rho_t(T_1, T_2) = \frac{\sum_{j} {\ell_j}  e_j(T_1-t) e_j(T_2-t) }{c_t(T_1) c_t(T_2) }, \qquad \mbox{for every }t< \tau,
	\end{equation*}
	with $c_t(T_i)$ as in Proposition \ref{prop:sde}.
\end{prop} 
\begin{proof}
	From equation \eqref{eq:Qexpansion}, accordingly to \cite[Prop. 4.3]{sdebook}, we derive the following expansion for the Wiener process $\W$:
	\begin{equation*}
		\W = \sum_j \sqrt{{\ell_j}} \beta^j e_j,
	\end{equation*}
	with $\beta^j := \frac{1}{\sqrt{{\ell_j}}}\langle \W, e_j \rangle_{}$ being real-valued mutually independent Wiener processes. Then, for every $t<\tau$ and every $x_1, x_2\ge0$, because $\delta_{x_i}$ is  a bounded and linear operator in $\h$, we obtain
		\begin{equation}
		\label{expansionWx}
		\delta_{x_i} \W_t = \W_t(x_i) = \sum_j \sqrt{{\ell_j}} \beta_t^j e_j(x_i) \in \R, \quad \mbox{for } i = 1, 2,
	\end{equation}
	which is a linear combination of real-valued Wiener processes, hence $\W_t(x_i)$ is a Gaussian random variable with mean $0$ and variance $\left(c_t(t+x_i)\right)^2t$. Then there exist two real-valued Wiener processes $W^1$ and $W^2$ such that 
	\begin{equation}
		\label{expansionWxuni}
		\W_t(x_i) = c_t(t+x_i) W^i_t, \quad \mbox{for } i = 1, 2.
	\end{equation}
	For $T_i = t+x_i$, we want now to find an expression for the correlation $\rho = \rho_t(T_1, T_2) := \mathrm{Corr}(W^1_t, W^2_t)$. By using the representation \eqref{expansionWx}, since the $\beta^j$'s are mutually independent, we first find that
	\begin{align}
		\notag
		\E\left[\W_t(T_1-t) \W_t(T_2-t)  \right] &= \E\left[\left( \sum_j \sqrt{{\ell_j}} \beta_t^j e_j(T_1-t)\right)\left( \sum_j \sqrt{{\ell_j}} \beta_t^j e_j(T_2-t)\right)  \right]\\
		&=    \sum_j {\ell_j}  e_j(T_1-t) e_j(T_2-t) \,t,\label{eq:cov1}
	\end{align}
	and, by using the representation \eqref{expansionWxuni} and that, by definition of $\rho$, we have $\E\left[ W^1_t W^2_t\right] = \rho_t(T_1, T_2) \,t$, we find that 
		\begin{equation}
		\label{eq:cov2}
		\E\left[\W_t(T_1-t) \W_t(T_2-t)  \right] = \E\left[ c_t(T_1) W^1_t c_t(T_2) W^2_t \right]=     c_t(T_1) c_t(T_2)\rho_t(T_1, T_2)  \,t.
	\end{equation}
	Comparing equations \eqref{eq:cov1} and \eqref{eq:cov2} we obtain the expression for the correlation $\rho_t(T_1, T_2)$.
\end{proof}

\subsection{Positivity of mild solutions}\label{sec:positivity}
It is required in some applications to work with strictly-positive curves. We derive in this section sufficient conditions for the SPDE coefficients \eqref{spde} to guarantee positivity of the solutions. More specifically, we study conditions guaranteeing solutions in $\h_+$ and conditions guaranteeing solutions in $\h_>$. We emphasize that guaranteeing solutions in $\h_>$ is crucial for the invertibility of the diffusion operator and hence, for example, for applying a change of measure from a physical to a risk-neutral measure. We also point out that positivity conditions for SPDEs  in the Filipovi\'c space have been studied before; for example \cite{filipovic2010} consider the case of bounded and Lipschitz coefficients. Generally speaking, it is usually required that the drift is pointing inwards at the boundary points of the cone $\h_+$, and that the diffusion is parallel to the boundary at the boundary points. Extensions of some results in \cite{filipovic2010} are in \cite{T2017} and \cite{T2022}. We provide here conditions which are tailor-made for coefficients defined based on point-wise operating maps, and which can be easily checked in this setting. 

Starting from the SPDE \eqref{spde}, we first consider the case with no drift, namely we set $\alpha =0$.
\begin{thrm}\label{prop:positivity_nodrift}
	Let $\psi\in\mathcal{C}^{1,2}(\R_{+}^0\times \R_{+}^0;\R_+^0)$
 satisfy the conditions of Lemma \ref{lipchitz:op:norm}. If:
	\begin{enumerate}
		\item For all $t\in\R^0_+$, $\psi(t, y)=0$ if and only if $y=0$;
		\item There exists an $\varepsilon>0$ such that for $\mathcal{O}_{\varepsilon}:=\{(t,y)| t\in\R^0 _+\mbox{ and } y\in (0,\varepsilon)\}$ it holds:
		\begin{enumerate}
			\item $\inf_{\mathcal{O}_{\varepsilon}}\frac{\partial_t\psi(t, y)}{\psi(t,y)}>-\infty$;
			\item $\inf_{\mathcal{O}_{\varepsilon}}\partial_y\psi(t, y)>-\infty$;
			\item $\inf_{\mathcal{O}_{\varepsilon}}\psi(t, y)\partial_{yy}\psi(t, y)>-\infty$.
		\end{enumerate}
	\end{enumerate} 
	Then for any initial condition $g_0\in \h_+$, respectively $g_0\in \h_>$, the (local) mild solution of the SPDE \eqref{spde} with multiplicative diffusion operator and zero drift belongs to $\h_+$, respectively to $\h_>$.
\end{thrm}
\begin{proof}
From Lemma \ref{lipchitz:op:norm} and \cite[Thm. 4.5]{T2012}, we know that there exists a unique (local) mild solution $(g, \tau)$ to \eqref{spde}, where $\tau$ is a stopping time almost surely positive and $g$ is given by
\begin{equation*}
	g_{t\land\tau} = \mathcal{S}_{t\land\tau} g _0 + \int_0^{t\land\tau}\mathcal{S}_{t\land\tau-s}\sigma (s, g_s) d\W_s,
\end{equation*}
where $\sigma (s, g_s) = \mathcal{M}_{\Psi (s, g_s)}$. In particular, $\tau = \infty$ if the solution is global. Moreover, from Proposition \ref{prop:sde} we know that for every $T\ge0$ there exists a one-dimensional Wiener process $W$ such that the forward price $F(\cdot, T)$ follows the dynamics \eqref{eq:sde_d} for any $t< \tau$. In particular, equation \eqref{eq:sde_d} is a one-dimensional SDE with strictly positive initial condition, $F(0, T)=g_0(T) >0$ for $g_0\in\h_+$. We shall now prove that $\{\{F(t\land \tau,T)\}_{t\in[0,T]}, \, T\ge 0\}$ is a strictly positive process. 

By assumption \emph{1.},  if $F(t\land \tau,T)=0$ for a certain $0\le t\le T$, then the process is absorbed in $0$, namely $F(u\land\tau, T) = 0$ for any $t\le u \le T$. Thus, to show that for all $T>0$, $F(t\land\tau,T)>0$ for all $0\le t\le T$ when $F(0,T)>0$, it suffices to show that $F(t\land\tau,T)$ cannot hit $0$ in finite time. For this, let us define the process $U(t,T) := \log( \psi(t, F(t,T)))$.  We then introduce $\tau_0^F:= \inf\{t\in\R_+^0 | \, F(t\land\tau,T)=0 \}$ and $\tau_{-\infty}^U:= \inf\{t\in \R_+^0 | \, U(t\land\tau,T)=-\infty \}$. We notice that $\tau_0^F<\infty$ 
if and only if $\tau_{-\infty}^U<\infty$, thus we shall prove that $\{\{U(t\land\tau,T)\}_{t\in[0,T]},\, T\ge 0\}$ does not explode in finite time when the value of $F(t\land\tau,T)$ is in a (right-)neighborhood of $0$. 

Let $\varepsilon>0$. Then, for $F(t,T)\in(0, \varepsilon)$ and $t<\tau$, we apply It\^o's formula to $U(t,T)$:
\begin{equation}\label{eq:dU}
	\begin{aligned}
	dU(t,T) =& \left(\frac{\partial_t\psi(t, F(t,T))}{\psi(t, F(t,T))}+\frac{c_t^2}{2}\left(\partial_{yy}\psi(t, F(t,T))\psi(t, F(t,T))-\left(\partial_y\psi(t, F(t,T))\right)^2\right)\right)dt \\&+ c_t\partial_y\psi(t, F(t,T)) dW_t.
	\end{aligned}
\end{equation}
By assumption \emph{2.}, the coefficients of the SDE \eqref{eq:dU} are bounded for $F(t,T)\in(0, \varepsilon)$. It then follows that $\abs{U(t,T)}$ is finite with probability 1 for $F(t,T)\in(0, \varepsilon)$. In particular it cannot reach $-\infty$ in finite time, namely $\tau_{-\infty}^U=\infty$ almost surely. This in turn implies that $\tau_0^F=\infty$, hence $F(t\land\tau,T)>0$ for all $0\le t\le T$ and every $T>0$. We conclude that $g_t\in\h_+$ for $g_0\in\h_+$. 

Let now $g_0\in\h_>$. Since $\h_> \subset \h_+$, we already know that, under Assumption \emph{1.} and \emph{2.}, $g_t\in\h_+$ for $g_0\in\h_>$. We then need to prove the additional property $\lim_{x\to \infty} g_t(x) > 0$ for all $t< \tau$, or, equivalently, that $\lim_{T\to \infty} F(t,T) > 0$ for all $t<\tau$. From the first part of the proof, we know that $F(t,T)$ can be written as in \eqref{eq:sde}, and that $U(t,T) = \log( \psi(t, F(t,T)))$ follows the dynamics \eqref{eq:dU}. This implies that $\lim_{T\to \infty} F(t,T) > 0$ if and only if $\lim_{T\to \infty} U(t,T) > -\infty$ for those curves $F(t, \cdot)$ such that $\lim_{T\to \infty} F(t,T)$ lies in a (right-)neighborhood of $0$. We shall then show this latter equivalent statement, namely for $\varepsilon>0$ and $\h^\varepsilon_>:= \{f\in\h_> | \lim_{x\to\infty}f(x) \in (0, \varepsilon)\}$, we shall prove that
\begin{equation}\label{eq:limU}
	\lim_{T\to \infty} U(t,T)  = \lim_{T\to \infty} \log( \psi(t, F(t,T)))> -\infty ,\qquad \mbox{for all } F(t,\cdot )\in\h_>^\varepsilon.
\end{equation}

From \eqref{eq:dU} we have that
\begin{equation*}\label{eq:U}
		U(t\land\tau,T) = U(0,T)+\int_0^{t\land\tau}\tilde{\mathrm{A}}(s,F(s,T))ds + \int_0^{t\land\tau}c_s\partial_y\psi(s, F(s,T)) dW_s,
\end{equation*}
with
\begin{align*}
	\tilde{\mathrm{A}}(s,F(s,T)) := \frac{\partial_t\psi(s, F(s,T))}{\psi(s, F(s,T))}+\frac{c_s^2}{2}\left(\partial_{yy}\psi(s, F(s,T))\psi(s, F(s,T))-\left(\partial_y\psi(s, F(s,T))\right)^2\right).
\end{align*}
Hence, in particular, \eqref{eq:limU} holds if and only if the following three conditions are satisfied:
\begin{enumerate}
	\item[(I)]  $\lim_{T\to \infty} U(0,T) > -\infty$ for all $F(t,\cdot )\in\h_>^\varepsilon$;
	\item[(II)]  $\lim_{T\to \infty} \int_{0}^{t\land\tau}\tilde{\mathrm{A}}(s,F(s,T))ds > -\infty$ for all $F(t,\cdot )\in\h_>^\varepsilon$;
	\item[(III)] $\lim_{T\to \infty} \int_0^{t\land\tau}c_s\partial_y\psi(s, F(s,T)) dW_s > -\infty$  for all $F(t,\cdot )\in\h_>^\varepsilon$.
\end{enumerate}
We notice that, under assumption \emph{1.}, condition (I) is trivially satisfied for $g_0\in\h_>$. We need to prove (II) and (III). We start with (II).

Let $\{T_n\}_{n}$ be a sequence such that $T_n\to \infty$ as $n\to\infty$. Then $\lim_{n\to\infty} F(t, T_n)\in(0, \varepsilon)$ for every $F(t,\cdot )\in\h_>^\varepsilon$, and $\lim_{T\to \infty} \int_{0}^{t\land\tau}\tilde{\mathrm{A}}(s,F(s,T))ds = \lim_{n\to \infty} \int_{0}^{t\land\tau}\tilde{\mathrm{A}}(s,F(s,T_n))ds$. Moreover, under assumption \emph{2.}, and because $0< c_s< \infty$ for all $s\ge 0$, there exists a uniform bound $M> 0$ such that
\begin{equation*}
|\tilde{\mathrm{A}}(s,F(s,T_n)) | \le M, \qquad \mbox{for all } s\ge 0 \mbox{ and for } n \mbox{ large enough.}
\end{equation*}
This allows us to apply the bounded convergence theorem to get that
\begin{equation*}
	\lim_{n\to \infty} \int_{0}^{t\land\tau}\tilde{\mathrm{A}}(s,F(s,T_n))ds =\int_{0}^{t\land\tau}\lim_{n\to \infty}\tilde{\mathrm{A}}(s,F(s,T_n))ds \ge -M\cdot (t\land\tau)  > -\infty,
\end{equation*}
Hence $\lim_{T\to \infty} \int_{0}^{t\land\tau}\tilde{\mathrm{A}}(s,F(s,T))ds > -\infty$ for all $F(t,\cdot )\in\h_>^\varepsilon$, and we conclude that (II) holds.

We are now left to prove (III). In particular, we notice that (III) holds if $$\lim_{T\to \infty} \E\left[\int_0^{t\land\tau}\left(c_s\partial_y\psi(s, F(s,T))\right)^2 ds\right] < \infty$$  for all $F(t,\cdot )\in\h_>^\varepsilon$, hence we shall show the latter instead. For this, we consider the same sequence $\{T_n\}_{n}$ as above. Then, similarly as before, under assumption \emph{2.}, and because $0< c_s< \infty$ for all $s\ge 0$, there exists a uniform bound $M> 0$ such that
	\begin{equation*}
		|c_s\partial_y\psi(s, F(s,T_n))| \le M, \qquad \mbox{for all } s\ge 0 \mbox{ and for } n \mbox{ large enough.}
	\end{equation*}
	This allows us to apply the bounded convergence theorem again to get that
	\begin{align*}
		\lim_{T\to \infty} \E\left[\int_0^{t\land\tau}\left(c_s\partial_y\psi(s, F(s,T))\right)^2 ds\right]  &= \lim_{n\to \infty}\E\left[\int_0^{t\land\tau}\left(c_s\partial_y\psi(s, F(s,T_n))\right)^2 ds\right]  \\&= \E\left[\int_0^{t\land\tau}\lim_{n\to \infty}\left(c_s\partial_y\psi(s, F(s,T_n))\right)^2 ds\right] \\&\le \E\left[\int_{0}^{t\land\tau}M^2ds\right] = M^2\cdot (t\land\tau)< \infty.
\end{align*}
We conclude that $g_t\in\h_>$ for $g_0\in\h_>$ and $t<\tau$. 
\end{proof}

	\begin{cor}\label{cor:positivity_betaphi}
		Let $\psi(t, y)  = \beta(t)\phi(y)$ with $\beta:\R_{+}^0\to \R_+$ differentiable and  let $\phi\in\mathcal{C}^{2}(\R_{+}^0;\R_+^0)$ 
  satisfy the conditions of Corollary \ref{lem:specialcase:mult:structure}. If:
	\begin{enumerate}
		\item $\phi(y)=0$ if and only if $y=0$;
		\item $\inf_{t\in\R_+^0}\frac{\beta'(t)}{\beta(t)}>-\infty$;
		\item There exists an $\varepsilon>0$ such that:
		\begin{enumerate}
			\item $\inf_{y\in (0,\varepsilon)}\phi'(y)>-\infty$;
			\item $\inf_{y\in (0,\varepsilon)}\phi(y)\phi''(y)>-\infty$;
		\end{enumerate}
	\end{enumerate} 
	then for any initial condition $g_0\in \h_+$, respectively $g_0\in \h_>$, the (local) mild solution of the SPDE \eqref{spde} with zero drift and with multiplicative diffusion operator $\sigma(t, f) = \mathcal{M}_{\Psi(t, f)}$, where $\Psi(t, f)(x) = \beta(t)\Phi(f)(x) =\beta(t) \phi(f(x))$, belongs to $\h_+$, respectively to $\h_>$, for every $t<\tau$.
\end{cor}
\begin{proof}
	We can easily check that for $\psi(t, y)  = \beta(t)\phi(y)$ and $\beta$ differentiable, conditions \emph{1.--3.} imply conditions \emph{1.} and \emph{2.} of Theorem \ref{prop:positivity_nodrift}.
\end{proof}

Theorem \ref{prop:positivity_nodrift} provides sufficient conditions to ensure that the solutions remain in $\h_+$ and $\h_>$ for an SPDE without drift. We now want to extend this result to SPDEs with drift. To do this, we shall introduce a measure change that transforms an SPDE without drift to an SPDE with generic drift $\alpha: \R_+^0 \times \h_+ \to \h$. We need some preliminaries. 

In the following, we say that a stochastic process $\varphi : \Omega \times \mathbb{R}_+^0\to \h$  satisfies the {\em generalized Novikov condition} if for $\bar{T}>0$ the following is fulfilled:
\begin{equation}\label{novikov:cond}
\mathbb{E} \left[\exp \left\{ \frac{1}{2} \int_0^{\bar{T}} \norm{ \varphi(s) }^2 ds\right\}\right] < \infty.
\end{equation}
Condition \eqref{novikov:cond} appears indeed as a generalization to Hilbert space-valued stochastic processes of the classical Novikov condition \cite[Thm. 8.6.5]{oksendal}. It is shown in \cite[Prop. 10.17]{sdebook} that \eqref{novikov:cond} is sufficient to  define a new measure $\tilde\p\sim \p$ on $(\Omega, \mathcal{F})$ with Radon–Nikodym derivative of the form
\begin{equation}
	\label{eq:radon}
\left.\frac{d \tilde\p}{d\mathbb{P}}\right|_{\F_{\bar{T}}}= \exp \left\{ \int_0^{\bar{T}} \langle \varphi(s), d \W_s \rangle - \frac{1}{2} \int_0^{\bar{T}} \norm{\varphi(s)}^2 ds \right\} .
\end{equation}
It follows then from \cite[Thm. 10.14 ]{sdebook}, that for every $0\le t\le \bar{T}$ the process $\tilde{\W}$ defined by 
\begin{equation}
	\label{eq:Wtilde}
\tilde{\W}_t := \W_t - \int_0^t \varphi(s) ds
\end{equation}
is a $Q$-Wiener process w.r.t. $\{ \mathcal{F}_t\}_{t\geq 0}$ on the probability space $(\Omega, \mathcal{F},\tilde\p)$.

For a given drift coefficient $\alpha$, the process $\varphi : \Omega \times \mathbb{R}_+^0\to \h$ that we need for our purposes is of the form 
\begin{equation}
\label{eq:sigmavarphi}
    \sigma(t,f)\varphi(t) = \alpha(t,f), \qquad \mbox{for } t\in\R_+^0, \,f\in \h.
\end{equation}
Since the volatility $\sigma$ is a multiplication operator, thanks to Proposition \ref{prop:invertibility} we can make $\varphi(t)$ explicit in equation \eqref{eq:sigmavarphi} if the operator $\sigma(t,f)$ is invertible. For this, we must restrict the kernel function of $\sigma(t,f)$ to be in $\h_< \cup \h_>$ for every $t \in \mathbb{R}_+^0$ and $f \in \h$.  In particular, this holds under the assumptions of Theorem \ref{prop:positivity_nodrift}, as we show in the next corollary.

\begin{cor}\label{cor:invertiblesigma2}
	Let $\psi:\R_{+}^0\times \R_{+}^0\to \R_+^0$ satisfy the conditions of Lemma \ref{lipchitz:op:norm} and of Theorem \ref{prop:positivity_nodrift}, and for any $g_0\in\h_>$ let $(g, \tau)$ denote the (local) mild solution of the SPDE with no drift. Then for every $t<\tau$ it holds $\Psi (t, g_t)\in \h_>$, and the diffusion operator $\sigma (t, g_t) = \mathcal{M}_{\Psi (t, g_t)}$ is invertible.
\end{cor}
\begin{proof}
	By definition  of the function $\psi$, we know that $\Psi (t, g_t)\in\h_+$  for every $t<\tau$. We need to show that $\lim_{x\to \infty}\Psi (t, g_t)(x) = \lim_{x\to \infty}\psi (t, g_t(x)) > 0$  for every $t<\tau$. However, from Theorem \ref{prop:positivity_nodrift} we know that $g_t\in\h_>$ , hence $\lim_{x\to \infty}g_t(x) > 0$ for every $t<\tau$. By continuity of $\psi$, we then get that
	\begin{equation*}
		\lim_{x\to \infty}\psi (t, g_t(x)) = \psi\left(t, \lim_{x\to \infty} g_t(x)\right) >0,
	\end{equation*}
	since, by assumption, $\psi(t, y)=0$ if and only if $y=0$. Hence $\Psi (t, g_t)\in \h_>$  for every $t<\tau$. Moreover, we know from Proposition \ref{prop:invertibility} that $\sigma (t, g_t) = \mathcal{M}_{\Psi (t, g_t)}$ is invertible if $\Psi (t, g_t)\in\h_>$. This concludes the proof.
\end{proof}

By Corollary \ref{cor:invertiblesigma2}, starting from \eqref{eq:sigmavarphi} we can now define explicitly $\varphi(t):=\sigma(t,f)^{-1} \alpha (t,f)$. Moreover, assuming that the Novikov condition \eqref{novikov:cond} holds, we can introduce by \eqref{eq:radon} the change of measure from $\p$ to an equivalent measure $\tilde\p$ in terms of the (local) mild solution $g$ of the SPDE without drift.  We use this change of measure in the following lemma to prove that we can guarantee solutions in $\h_+$ and in $\h_>$ also for the SPDE with non-zero drift.

\begin{lem}\label{lem:positivity}
	Let $\psi:\R_{+}^0\times \R_{+}^0\to \R_+^0$ satisfy the assumptions of Theorem \ref{prop:positivity_nodrift}, and, for any initial condition $g_0\in \h_>$, let $(g, \tau)$ denote the (local) mild solution from Theorem \ref{prop:positivity_nodrift} of the SPDE without drift.  Let further $\alpha:\R_+^0\times \h_+\to \h$ satisfy either locally Lipschitz and locally bounded conditions, or locally Lipschitz and linear growth conditions, accordingly to the results in Section \ref{sec:conditions}. For any initial condition $\tilde g_0\in \h_>$, let then $(\tilde g , \tilde\tau)$  denote the (local) mild solution of the SPDE with drift $\alpha$ with $\tilde\tau\le \tau$ almost surely. If the condition 
	\begin{equation}\label{novikov:cond_positivity}
		\mathbb{E} \left[\exp \left\{ \frac{1}{2} \int_0^{\bar{T}} \norm{ \sigma(t,g_t)^{-1}\alpha(t, g_t) }^2 dt\right\}\right] < \infty
	\end{equation}
	is satisfied for some $\bar{T}<\tilde\tau$, then $\tilde g_{t}\in \h_>$ for any $t<\bar T$.
\end{lem}
\begin{proof}
	From Theorem \ref{prop:positivity_nodrift} we know that for any $g_0\in\h_>$ the (local) mild solution of
	\begin{equation*}
		dg_t = \partial_xg_tdt + \sigma(t, g_t)d\W_t
	\end{equation*}
	belongs to $\h_>$ for every $t< \tau$. Moreover, from Corollary \ref{cor:invertiblesigma2} we know that the diffusion operator $\sigma (t, g_t)$ is invertible for every $t< \tau$. For a given drift coefficient $\alpha$, we can then introduce the process $\varphi : \Omega \times \mathbb{R}_+^0\to \h$ by
	\begin{equation*}
		\varphi(t) :=\sigma(t,f)^{-1}\alpha(t, g_t), \quad \mbox{for all } t< \tau.
	\end{equation*}
    Assumption \eqref{novikov:cond_positivity} is the generalized Novikov condition \eqref{novikov:cond} for $\varphi$ defined above. We can then introduce by means of $\varphi$ a new measure $\tilde\p\sim \p$ and a Brownian motion $\{\tilde{\W}_t\}_{t\geq 0}$ as in \eqref{eq:radon} and \eqref{eq:Wtilde}, respectively. The process $g_t$ has the following dynamics under $\tilde\p$:
	\begin{equation}
		\label{spde:change3}
		d g_t = \partial _x  g_t dt +\alpha(t,  g_t)dt +\sigma (t, g_t) d \tilde{\W}_t.
	\end{equation}
    By the assumptions on $\alpha$ and $\sigma$, we know from \cite[Thm. 4.5]{T2012} that there exists a unique (local) mild solution to \eqref{spde:change3}. That means that $\tilde g$ must coincide with $g$. 
    Since the solution $\{g_t\}_{t\ge 0}$ of \eqref{spde:change3} is in $\h_>$ for any initial condition $g_0\in\h_>$, we conclude that also $\{\tilde g_t\}_{t\ge 0}$ is in $\h_>$ for any initial condition $\tilde g_0\in\h_>.$
\end{proof}

\begin{remark}
	We stress that, thanks to Lemma \ref{lem:positivity}, we can define a process $\varphi$ to change measure from $\p$ to an equivalent measure $\tilde\p$ in terms of the (local) mild solution of the SPDE with drift. In this way, we can change to an equivalent martingale measure, which is important for pricing derivatives when the initial model is set up under the physical measure $\p$.
\end{remark}

\section{Specifications}
\label{sec:spec}
In this section we consider some specifications for the SPDE coefficients \eqref{spde}, and we apply the results from the previous sections to show existence and uniqueness of (local) mild solutions, as well as positivity of the solutions. We consider first a class of exponential models, and, by deriving their dynamics, we reveal that they can be seen as a Hilbert-space valued counterpart of one-dimensional geometric models, as for example the Black-Scholes model. The second class of models we consider is a Hilbert-space valued counterpart of a constant elasticity of variance (CEV) model.

\subsection{Exponential models}
\label{sec:exponential}
We consider a class of exponential models. Let us start with the following two lemmas, where $\exp:\R\to\R_{+}$ and $\log:\R_{+}\to\R$ denote, respectively, the exponential and the logarithm maps on real numbers.
\begin{lem} 
	\label{lem:exp}
	The following statements are true:
	\begin{enumerate}
		\item The exponential map $\Exp: f \mapsto  \{ x\mapsto \exp(f(x))  \}$ is well defined as a map from $\h$ to $\h_>$.
		\item The logarithm map $\Log: h \mapsto  \{ x\mapsto \log(h(x))  \}$ is well defined as a map from $\h_>$ to $\h$.
	\end{enumerate}
Moreover, the map $\Exp$ acts as a left inverse to $\Log$ on $\h_>$, namely
$\Exp (\Log (h)) = h$ for all $h \in \h_>$, and the map $\Log$ acts as a left inverse to $\Exp$ on $\h$, namely $\Log (\Exp (f)) = f$ for all $f \in \h.$
\end{lem}
\begin{proof}
It is shown in \cite[Lemma 3.12]{fredkru2} that $\Exp(f)$ is well defined as an element in $\h$. The same follows also from  Proposition~\ref{prop:mapping:into:h}. The fact that $\Exp(f) \in \h_>$ is an easy consequence of Lemma~\ref{bounded:fcts:in:h}, the convergence $\lim_{x\rightarrow \infty } f(x)$, and the continuity of the exponential function. 

For \emph{2.}, since $h(x)>0$ for all $h\in \h_>$, we notice that $\log (h(x))$ is well defined for every $x\geq 0$. Moreover, there exists $\varepsilon>0$ such that  $h(x)\geq \varepsilon$ for every $x\geq 0$. Then there exists also a constant $C>0$ such that $1/h(x)<C$ for all $x$ and $|\log(h(0))| < C$. Combining, we get: 
\begin{align}	
	\norm{\Log(h)}^2_{}& =\log (h (0))^2 + \int_0^{\infty}  \left(\left(\log(h(x))\right)'\right)^2 w (x)dx \nonumber  \\
	& = \log (h (0))^2 + \int_0^{\infty} \left(\frac{h'(x)}{h(x)}\right)^2  w (x)dx \nonumber \\
		& \leq  C^2 + C^2 \int_0^{\infty}  \left(h'(x)\right)^2 w (x)dx  \leq  C^2\left(1+ \norm{h}^2\right) <  \infty, \nonumber
\end{align}	
	from which it follows that $\Log(h)\in \h$ and the map $\Log : \h_> \rightarrow \h$ is therefore well defined. 

The last statement follows from the real-valued case, because the two operators act point-wise.
\end{proof}

Let now $(g, \tau)$ be a (local) mild solution of the SPDE \eqref{spde} for an initial condition $g_0\in \h$, 
and let $z_t := \Exp(g_t)$. Then $\{z_t\}_{t<\tau}\in\h_>$ by Lemma \ref{lem:exp}. Moreover, we can prove the following result for the dynamics of $z_t$.
\begin{thrm}\label{thrm:geometric}
Let $(g, \tau)$ be a (local) mild solution to the SPDE \eqref{spde}, and let $z_t := \Exp(g_t)$ for all $t<\tau$. Then the process $\{z_t\}_{t<\tau}$ solves the following SPDE: 
	\begin{equation*}
		dz_t = \partial _x z_t dt +  z_t\,\tilde\alpha(t, z_t) dt + z_t\,\tilde\sigma (t, z_t) d\W_t,
	\end{equation*} 
	with $\tilde\alpha: \R_{+}^0\times \h_> \rightarrow \h$ and $\tilde\sigma: \R_{+}^0\times \h_> \rightarrow \mathcal{L}(\h)$ given by 
	\begin{align*}
		\tilde\alpha(t, h) := \alpha(t, \Log(h))+\frac{1}{2}\Sigma(t, \Log(h)),\quad \mbox{ and } \quad \tilde\sigma(t, h) := \sigma(t, \Log(h)),
	\end{align*} with
		\begin{equation}
		\label{Sigmahat}
		\Sigma(t, f):= \sigma(t, f)\q\,\sigma(t, f)^* \delta_{T-t}^*(1), 
	\end{equation}
	for all $f\in \h$ and $t<\tau$, and  $\alpha$ and $\sigma$ as in equation \eqref{spde}.
\end{thrm}
\begin{proof}
First note that because $g_t$ is not an It\^o process, we can not directly apply It\^o's formula to $\Exp(g_t)$ to compute the dynamics of $z_t$. Hence we first change the parametrization, namely we consider $f(t, T)= g_t(T-t)$. Then, we transform the stochastic integral driven by the $\h$-valued Wiener process $\W$ into a stochastic integral driven by a one-dimensional Wiener process. At this point, we are dealing with a diffusion process and we can apply It\^o's formula. After that, we reverse the two preliminary steps to obtain our result.

Let $g_t$ be the mild solution in \eqref{mildsol}
for an initial condition $g_0\in \h$. Then the forward price $f(t, T)$ for $T\ge t$ is given by $f(t, T)= g_t(T-t) = \delta_{T-t}g_t$. In particular, since the point-wise evaluation $\delta_x$ is a continuous linear functional on $\h$ for all $x\in \R_{+}^0$, see \cite[Thm. 5.1]{filipovic2010}, we can interchange integration and point-wise evaluation. We then write:
\begin{align}
	f(t, T) &=\delta_{T-t}\mathcal{S}_t g_0+ \int_0^t \delta_{T-t}\semi_{t-s}\alpha (s, g_s) ds + \int_{0}^{t}\delta_{T-t}\semi_{t-s}\sigma(s, g_s) d\W_s\notag\\
	&=\delta_{T}g_0 + \int_0^t \delta_{T-s}\alpha(s, g_s) ds + \int_{0}^{t}\delta_{T-s}\sigma(s, g_s) d\W_s\notag\\
	&=\delta_{T}g_0 + \int_0^t \delta_{T-s}\alpha(s, f(s, s+\cdot)) ds + \int_{0}^{t}\delta_{T-s}\sigma(s, f(s, s+\cdot)) d\W_s,\label{eq:ftT}
\end{align}
where we used that $\delta_x\semi_s = \delta_{x+s}$ and $\delta_x g_s= g_s(x) = f(s, s+x) = \delta_xf(s, s+\cdot)$.  

By proceeding as in the proof of Proposition \ref{prop:sde}, we introduce a one-dimensional Wiener process $W$ such that 
\begin{equation}
	\label{eq:transformation}
	\delta_{T-s}\sigma(s, f(s, s+\cdot)) d\W_s = \left( \delta_{T-s}\sigma(s, f(s, s+\cdot))  \q \sigma(s, f(s, s+\cdot))^* \delta_{T-s}^*\left(1\right)\right)^{1/2} dW_s,
\end{equation}
and for every $0\le t \le T$, we rewrite the last term in equation \eqref{eq:ftT} as
\begin{equation*}
	\int_{0}^{t}\delta_{T-s}\sigma(s, f(s, s+\cdot)) d\W_s =\int_{0}^{t}\Gamma(s, f(s, s+\cdot))dW_s,\\
\end{equation*}
where 
\begin{equation}
	\label{sigma}
	\Gamma(s, h)^2:= \delta_{T-s}\sigma(s, h)\q\,\sigma(s, h)^* \delta_{T-s}^*(1), \qquad \mbox{for every }h\in \h.
\end{equation}
Equation \eqref{eq:ftT} now reads like
\begin{equation*}
		f(t, T) =\delta_{T}g_0 + \int_0^t \delta_{T-s}\alpha(s, f(s, s+\cdot)) ds + \int_{0}^{t}\Gamma(s, f(s, s+\cdot))dW_s.
\end{equation*}
We then apply It\^o's formula to $y_t(T):=\exp(f(t, T))$. This gives us
\begin{equation*}
\begin{aligned}
	dy_t(T) =&\, y_t(T) \biggl( \Bigl( \delta_{T-t}\alpha(t, \Log \left(y_t(t+\cdot)\right))+ \frac{1}{2}\Gamma(t, \Log \left(y_t(t+\cdot)\right)) ^2\Bigr)dt\biggr.\\&  \biggl.+ \Gamma(t, \Log \left(y_t(t+\cdot)\right)) dW_t \biggr),
\end{aligned}
\end{equation*}
where $f(t, t+\cdot) = \Log \left(y_t(t+\cdot)\right)$ is well defined because of Lemma \ref{lem:exp}. By rewriting the last equation in integral form and substituting $T=t+x$, we obtain
\begin{align*}
		y_t(t+x) =&\, y_0(t+x)  + \int_0^t y_s(t+x) \Bigl( \delta_{t+x-s}\alpha(s, \Log \left(y_s(s+\cdot)\right)) \Bigr.\\ 
        &\Bigl. +\frac{1}{2}\delta_{t+x-s}\Sigma(s, \Log \left(y_s(s+\cdot)\right))\Bigr)ds \\ 
        &\!+ \int_0^t y_s(t+x)\Gamma(s, \Log \left(y_s(s+\cdot)\right)) dW_s,
\end{align*}
with $\Sigma$ in equation \eqref{Sigmahat} such that $\delta_{T-s}\Sigma(s, h)= \Gamma(s, h)^2$ for every $h\in\h$.

By definition of the process $z$, it holds that $z_t(x) = y_t(t+x)$. We can then write
\begin{equation}
	\label{eq:ztx}
	\begin{aligned}
	z_t(x) =&\,  z_0(t+x)  + \int_0^t \delta_{t+x-s}z_s \Bigl( \alpha(s, \Log\left(z_s\right)) +\frac{1}{2}\Sigma(s, \Log\left(z_s\right)) \Bigr)ds\\& + \int_0^t\left(\delta_{t+x-s}z_s\right)\Gamma(s, \Log\left(z_s\right)) dW_s.
	\end{aligned}
\end{equation}
We now substitute \eqref{sigma} with $T=t+x$ into the stochastic integral in equation \eqref{eq:ztx}:
\begin{align}
	&\int_0^t\left(\delta_{t+x-s}z_s\right)\Gamma(s, \Log\left(z_s\right)) dW_s \notag\\&= \int_0^t\left(\delta_{t+x-s}z_s\right)\left( \delta_{t+x-s}\sigma(s, \Log\left(z_s\right))\q\,\sigma(s, \Log\left(z_s\right))^* \delta_{t+x-s}^*(1)\right)^{1/2} dW_s\label{eq:proof1}\\
	&= \int_0^t \left(\delta_{t+x-s}z_s\right)\delta_{t+x-s}\sigma(s,\Log\left(z_s\right))d\W_s = \int_0^t \delta_{t+x-s}z_s\sigma(s, \Log\left(z_s\right))d\W_s,\label{eq:proof2}
\end{align}
where from \eqref{eq:proof1} to  \eqref{eq:proof2} we used \eqref{eq:transformation}. Substituting then \eqref{eq:proof2} into  \eqref{eq:ztx} we obtain that
\begin{align*}
	z_t = \semi_t z_0 +\int_0^t \semi_{t-s}z_s \bigl(  \alpha(s, \Log\left(z_s\right))+ \frac{1}{2}\Sigma(s, \Log\left(z_s\right)) \bigr)ds+\int_0^t \semi_{t-s}z_s\sigma(s, \Log\left(z_s\right))d\W_s,
\end{align*}
hence $z_t$ satisfies the SPDE
\begin{equation*}
	dz_t = \partial _x z_t dt +   z_t\bigl(\alpha\left(t, \Log\left(z_t\right)\right)+\frac{1}{2}\Sigma\left(t, \Log\left(z_t\right)\right) \bigr) dt + z_t \sigma \left(t, \Log\left(z_t\right)\right)d\W_t,
\end{equation*}
which concludes the proof.
\end{proof}

\begin{remark}
	If the coefficients in \eqref{spde} are deterministic, then the process $z_t$ follows the following dynamics:
		\begin{equation*}
		dz_t = \partial _x z_t dt +  z_t \Bigl(\alpha(t)+\frac{1}{2}\sigma(t)\q\,\sigma(t)^* \delta_{T-t}^*(1)\Bigr)dt + z_t  \sigma(t) d\W_t.
	\end{equation*} 
	In this setting, one can, for example, compute explicit price formulas for European options on energy, see \cite[Section 3.3]{fredkru2}.
\end{remark}

From Theorem \ref{thrm:geometric}, we see that an exponential model solves an SPDE of the form
\begin{equation}
\label{spdez:general}
dz_t = \partial _x z_t dt +  z_t\,\tilde\alpha(t, z_t) dt + z_t\,\tilde\sigma (t, z_t) d\W_t,
\end{equation} 
for some $\tilde\alpha: \R_{+}^0\times \h_> \rightarrow \h$ and $\tilde\sigma: \R_{+}^0\times \h_> \rightarrow \mathcal{L}(\h)$. The structure of the SPDE \eqref{spdez:general} resembles that of a geometric model in finite dimension. In particular, equation \eqref{spdez:general} corresponds to the SPDE \eqref{spde} with coefficients of the form $\alpha(t,h)= h\, \tilde{\alpha}(t,h)$ and $\sigma(t,h)= h\, \tilde{\sigma}(t,h)$ for all $h\in\h_>$ and $t\in\R^0_+$. Then, if $\sigma(t,h) = \mathcal{M}_{\Psi(t,h)}$ with $\Psi(t,h)(x)= \psi(t,h(x))$ for some $\psi: \R_{+}^0\times\R_{+} \rightarrow  \mathbb{R}$, we can introduce $\tilde\psi: \R_{+}^0\times\R_{+} \rightarrow  \mathbb{R}$ such that $\psi(t, y)= y\,\tilde{\psi}(t, y)$ 
and $\tilde\sigma(t,h) = \mathcal{M}_{\tilde\Psi(t,h)}$ for $\tilde\Psi(t,h)(x)= \tilde\psi(t,h(x))$ 
for all $h\in\h_>$ and $t\in\R^0_+$. Similarly, we can introduce $\tilde a: \R_{+}^0\times\R_{+} \rightarrow  \mathbb{R}$ such that $\tilde\alpha(t,h)(x) = \tilde a(t,h(x))$ and $\alpha(t,h)(x) = h(x)\,\tilde a(t,h(x))$.

\begin{cor}
\label{cor:exponentialmodel}
Let $\psi: \R_{+}^0\times\R_{+} \rightarrow  \mathbb{R}$ and $\tilde\psi: \R_{+}^0\times\R_{+} \rightarrow  \mathbb{R}$ be such that $\psi(t, y)= y\,\tilde{\psi}(t, y)$ for all $t\in\R_+^0$ and $y\in\R_+$. Then the following statements hold:
\begin{enumerate}
    \item Let $\tilde\psi $ be differentiable in the second variable and let the functions $\tilde\psi$ and $\partial_y\tilde\psi $ be continuous in $\R_+^0 \times \R_+$. Let further $\partial_y \tilde\psi: \mathbb{R}_+^0 \times \mathbb{R}_+ \rightarrow  \mathbb{R}$ be locally Lipschitz. Then $(t,f) \mapsto \Psi(t, f)$ and $(t,f)\mapsto  \mathcal{M}_{\Psi (t, f)}$ are locally bounded and locally Lipschitz.
    
    \item Let $\sup_{(t,y) \in  \R_{+}^0\times\R_{+} }  \tilde{\psi}(t, y)+y\,\partial_y \tilde\psi(t,y)<\infty$. Then $(t,f) \mapsto \Psi(t, f)$ and $(t,f)\mapsto  \mathcal{M}_{\Psi (t, f)}$ are of linear growth.
\end{enumerate}
\end{cor}
\begin{proof}
The corollary is an application of Proposition \ref{prop:lipchitz}, Proposition \ref{prop:bounded} and Proposition \ref{prop:lin:growth} to maps $\psi$ of the form $\psi(t, y)= y\,\tilde{\psi}(t, y)$ with $\partial_y\psi(t, y)= \tilde{\psi}(t, y)+y\,\partial_y \tilde\psi(t,y)$.
\end{proof}

In particular, when the process $\{z_t\}_{t<\tau}$ is defined by $z_t = \Exp(g_t)$ starting from $(g, \tau)$ solution to the SPDE \eqref{spde}, we know from Lemma \ref{lem:exp} that $\{z_t\}_{t<\tau}\in\h_>$. Let now the process $(z, \tau)$ be a solution to \eqref{spdez:general} corresponding to some initial condition $z_0$. 
Under suitable conditions on $\tilde\sigma$ and $\tilde\alpha$, we obtain again that $\{z_t\}_{t<\tau}\in\h_>$ for all $z_0\in \h_>$.  

\begin{cor}\label{cor:posExp}
	Let  $\tilde\psi\in\mathcal{C}^{1,2}(\R_{+}^0\times \R_{+};\R_+)$ 
 satisfy the conditions of Corollary \ref{cor:exponentialmodel}. If there exists an $\varepsilon>0$ such that for $\mathcal{O}_{\varepsilon}:=\{(t,y)| t\in\R^0 _+\mbox{ and } y\in (0,\varepsilon)\}$ it holds:
		\begin{enumerate}
			\item $\inf_{\mathcal{O}_{\varepsilon}}\frac{\partial_t\tilde\psi(t, y)}{\tilde\psi(t,y)}>-\infty$;
			\item $\inf_{\mathcal{O}_{\varepsilon}}y\,\partial_y\tilde\psi(t, y)>-\infty$;
			\item $\inf_{\mathcal{O}_{\varepsilon}}y\,\tilde\psi(t, y)\partial_{yy}\tilde\psi(t, y)>-\infty$.
		\end{enumerate} 
	Then for any initial condition $z_0\in \h_>$, the (local) mild solution of the SPDE \eqref{spdez:general} with multiplicative diffusion operator and zero drift belongs to $\h_>$.

    Let further $\tilde a:\R_+^0\times \R_+\to \R$ satisfy the conditions of Corollary \ref{cor:exponentialmodel}, and for any initial condition $\bar z_0\in \h_>$ denote by $(\bar z , \bar\tau)$ the (local) mild solution of the SPDE \eqref{spdez:general} with drift $\tilde\alpha(t,h)(x) = \tilde a(t,h(x))$ where $\bar\tau\le \tau$ almost surely. If the condition 
	\begin{equation*}
		\mathbb{E} \left[\exp \left\{ \frac{1}{2} \int_0^{\bar{T}} \norm{ \tilde\sigma(t, z_t)^{-1}\tilde\alpha(t, z_t) }^2 dt\right\}\right] < \infty
	\end{equation*}
	is satisfied for some $\bar{T}<\tau$, then $\bar z_{t}\in \h_>$ for any $t<\bar T$.
\end{cor}
\begin{proof}
    This is an application of Theorem \ref{prop:positivity_nodrift} for the first part of the statement, and of Lemma \ref{lem:positivity} for the second part, to maps $\psi(t,y) = y\,\tilde\psi(t,y)$, $\sigma(t,h) = h\tilde\sigma(t,h)$ and $\alpha(t,h) = h\tilde\alpha(t,h)$.
\end{proof}

Thanks to Corollary \ref{cor:posExp} we know that, for any initial condition $z_0\in \h_>$, the solution $(z, \tau)$ of the SPDE \eqref{spdez:general} is guaranteed to be in $\h_>$ for suitable $\tilde\alpha: \R_{+}^0\times \h_> \rightarrow \h$ and $\tilde\sigma: \R_{+}^0\times \h_> \rightarrow \mathcal{L}(\h)$. Moreover, from Lemma \ref{lem:exp}, the map $\textrm{Log}$ is well defined on $\h_>$. This allows us to define, starting from $\{z_t\}_{t<\tau}$, the process $g_t:=\Log(z_t)$. Furthermore, one can derive the dynamics of $\{g_t\}_{t<\tau}$ with similar steps to the ones in the proof of Theorem \ref{thrm:geometric}.

\subsection{The CEV specification}
\label{sec:cev}
We consider a constant elasticity of variance (CEV) specification for the volatility operator. In particular, for $\gamma\ge0$ we consider a model of the form \eqref{eq:betaphi} with $$\Phi(f) = \Phi^\gamma(f):= f^\gamma= \left\{x\mapsto \phi^\gamma(f(x)):=f(x)^{\gamma}\right\},$$ where (for now) $\phi^\gamma:\R\to\R$ and $\Phi^\gamma:\h\to\h$. We then define $\psi^\gamma(t, y) := \beta(t)\phi^\gamma(y) = \beta(t)y^\gamma$, for $t\in\R_{+}^0$ and $y\in\R$. Here $\beta:\R_+^0\to\R_+$ is a deterministic, continuous and bounded function modeling, e.g., the  seasonality in the volatility.  
The kernel function for the multiplicative operator is of the form
\begin{equation*}
	\begin{aligned}
		\Psi^\gamma: \R_{+}^0\times \h & \to \h \\
(t, f)&\mapsto  \{x\mapsto \beta(t) f(x)^\gamma\},
	\end{aligned}
\end{equation*}
and the volatility operator is then given by
\begin{equation*}
	\sigma_{\mathrm{cev}}(t, f)  := \mathcal{M}_{\Psi^\gamma(t, f)} =\mathcal{M}_{\beta(t)\Phi^{\gamma}(f)}= \mathcal{M}_{\beta(t)f^{\gamma}}.
\end{equation*}
We shall identify the (multiplicative) volatility operator with its kernel function. With the CEV specification, the SPDE \eqref{spde} becomes
\begin{equation}
	\label{spde_cev}
	dg_t = \partial _x g_t dt + \alpha(t, g_t)dt + \beta(t)g_t^{\gamma} d\W_t.
\end{equation} 
We refer to equation \eqref{spde_cev} as the CEV-SPDE.

As in the standard CEV model, the parameter $\gamma$ controls the relationship between volatility and curve level. 
For $\gamma = 1$ we obtain a model corresponding to an infinite-dimensional counterpart of the geometric Brownian motion, as discussed, e.g., in \cite[Example 2.7]{barth14} and covered in Section \ref{sec:exponential} when $ \tilde\sigma(t, z_t) = \beta(t)$ in equation \eqref{spdez:general}. For $\gamma = 0$ one obtains a deterministic volatility operator.

\begin{remark}
	It is also possible to generalize the CEV model \eqref{spde_cev} by considering real-valued time-dependent and state-dependent exponents of the form $\gamma = \tilde\gamma(t,f(\cdot))$ for $t\in\R_{+}^0$, $f\in\h$ and $\tilde\gamma:\R_+^0\times\R\to\R_{+}^0$. However, for the present study we restrict ourselves to deterministic and constant $\gamma$. An exception is done in Proposition \ref{prop:gamma12} in order to ensure Lipschitz continuity for $1<\gamma<2$.
\end{remark}

For $\sigma_{\mathrm{cev}}$ to be well defined, we need $\mathcal{M}_{\Psi^\gamma(t, f)}\in\mathcal{L}(\h)$ for all $t\in\R_+^0$ and $f\in\h$. Moreover, we need $(t,f)\mapsto\mathcal{M}_{\Psi^\gamma(t, f)}$ to be locally Lipschitz and of linear growth (locally bounded) in order to ensure existence and uniqueness of (local) mild solutions  for the CEV-SPDE. Hence we must derive conditions on the exponent $\gamma$ in order to guarantee these properties for $\mathcal{M}_{\Psi^\gamma(t, f)}$. However, from Remark \ref{rmk:phi}, we deduce that if $\mathcal{M}_{f^{\gamma}}\in \mathcal{L}(\h)$, then $\mathcal{M}_{\Psi^\gamma(t, f)}\in\mathcal{L}(\h)$ for every $f\in\h$ and $t\in\R_+^0$. Moreover,  if $\mathcal{M}_{f^{\gamma}}$ is locally Lipschitz and/or locally bounded and/or of linear growth, then also $\mathcal{M}_{\Psi^\gamma(t, f)}$ is locally Lipschitz and/or locally bounded and/or of linear growth for every $f\in\h$ and $t\in\R_+^0$. For the rest of the section, we shall therefore focus on studying the map $\Phi^\gamma(f) = f^\gamma$ and the multiplicative operator $\mathcal{M}_{f^{\gamma}}$.

The fact that $f\mapsto\mathcal{M}_{f^{\gamma}}$ is well defined in $\mathcal{L}(\h)$ and locally bounded follow from Proposition \ref{prop:mapping:into:h} and Proposition \ref{prop:locallybounded}.
	\begin{cor}
		\label{lem:fgamma2}
		The following statements hold: 
		\begin{enumerate}
			\item If $\gamma \in \mathbb{N}$, then $\Phi^\gamma(f)=f^{\gamma} \in \h$ for all $f\in \h$, and $\Phi^\gamma$ is locally bounded in $\h$. Moreover, the operator $\mathcal{M}_{_{f^{\gamma}}}$ is well defined and the map $f\mapsto \mathcal{M}_{f^{\gamma}}$ is locally bounded in $\h$.
		\item If $\gamma \in \mathbb{R}$ with $\gamma \geq 1$, then $\Phi^\gamma(f)=f^{\gamma} \in \h^0_+$ for all $f\in \h_+^0$, respectively $\Phi^\gamma(f)=f^{\gamma} \in \h_+$ for all $f\in\h_+$, and $\Phi^\gamma$ is locally bounded. Moreover, the operator $\mathcal{M}_{_{f^{\gamma}}}$ is well defined and the map $f\mapsto \mathcal{M}_{f^{\gamma}}$ is locally bounded both in $\h_+^0$ and in $\h_+$.
		\end{enumerate}
	\end{cor}
	\begin{proof}
		For $\gamma \in \mathbb{N}$,  the function $\phi^\gamma: \mathbb{R} \rightarrow \mathbb{R}$ defined by $\phi^\gamma(y)=y^{\gamma}$ is differentiable in $\R$, hence satisfies case \emph{1.} of Proposition \ref{prop:mapping:into:h} and case \emph{1.} of Proposition \ref{prop:locallybounded}. Then  $f^{\gamma} \in \h$ for $f\in \h$ and the map $f\mapsto f^\gamma$ is locally bounded in $\h$. By \eqref{norm:mult:operator}, the map $f\mapsto \mathcal{M}_{f^{\gamma}}$ is also well defined and locally bounded in $\h$. 
		
        Let now $\mathcal{X}\in \{\R_+^0, \R_+\}$ and $\mathcal{Y}\in\{\h_+^0, \h_+\}$. Then, for $\gamma \in \mathbb{R}$ with $\gamma\geq 1$, the function $\phi^\gamma: \mathcal{X} \rightarrow \mathcal{X}$ defined by $\phi^\gamma(y)=y^{\gamma}$ fulfills the assumptions of case \emph{2.} of Proposition \ref{prop:mapping:into:h} and of cases \emph{2.} and \emph{3.} of Proposition \ref{prop:locallybounded}, hence  $f^{\gamma} \in \mathcal{Y}$ for $f\in \mathcal{Y}$ and the map $f\mapsto f^\gamma$ is locally bounded in $\mathcal{Y}$. By \eqref{norm:mult:operator}, the map $f\mapsto \mathcal{M}_{f^{\gamma}}$ is also well defined and locally bounded in $\mathcal{Y}$. 
	\end{proof}

From Corollary \ref{lem:fgamma2}, we get that the multiplicative operator $\mathcal{M}_{_{f^{\gamma}}}$ is well defined and locally bounded in $\h_+^0$ and $\h_+$ only for $\gamma \ge 1$. Indeed, we notice that for $0< \gamma<1$ the condition $\lim_{y\to0} \left(\phi^{\gamma}\right)^\prime(y) = \gamma\lim  y^{\gamma-1}<\infty$ of case \emph{2.} of Proposition \ref{prop:mapping:into:h} and of cases \emph{2.} and \emph{3.} of Proposition \ref{prop:locallybounded} is not satisfied. Hence we need to restrict ourselves to $\gamma \ge1$. The main reason for this lies in the definition of the norm in the Filipovi\'c space which involves the derivative, namely the norm of $h\in\h$ involves $h'$, hence the norm of $\Phi^\gamma(f) = \phi^\gamma(f(\cdot))$ involves the derivative $\left(\phi^{\gamma}\right)^\prime$.  Indeed, we know that $f\in \h$ is weakly differentiable by definition. It follows that also $f^{\gamma}$ is weakly differentiable, and by the chain rule we obtain that
	\begin{equation}
		\label{normfgamma}
		\norm{f^{\gamma}}^2 = \left(f(0)^{\gamma} \right)^2  + \int_0^{\infty} \left(\gamma  f(x)^{\gamma-1} f'(x)\right)^2 w (x)dx.
	\end{equation} 
	Here we can see that for $\gamma < 1$ the term $f(x)^{\gamma-1}$ gets large for $x$ such that $f(x)$ is small. As a result, the norm $\norm{f^{\gamma}}_{}$ might $\infty$. Hence we need to restrict ourselves to $\gamma \ge 1$. 

We further notice that, as in the one-dimensional case, the linear growth condition is not satisfied by $\phi^\gamma$ for any $\gamma > 1$. Hence, accordingly to Corollary \ref{lem:specialcase:mult:structure3}, the map $f\mapsto\mathcal{M}_{f^\gamma}$ is of linear growth only for $\gamma = 1$, and we can get existence and uniqueness of global mild solutions only in this case.

\begin{remark}\label{remark:noboundf}
We point out that for $f\in \h_> $, we know from Lemma~\ref{bounded:fcts:in:h} that $\inf_{x\in \mathbb{R}_+^0}\abs{f(x)} > 0$ and thus $\sup_{x\in \mathbb{R}_+^0}\abs{f(x)^{\gamma-1}} < \infty$. This means that for $f\in \h_>$ the map $f\mapsto\mathcal{M}_{f^\gamma}$ is well defined also for $0<\gamma < 1$, as we can deduce from equation \eqref{normfgamma}. However, as pointed out in Remark \ref{remark:infh}, we do not have a lower bound for $\inf_{x\in \mathbb{R}_+^0}\abs{f(x)}$ in terms of the norm of $f$, hence we do not have an upper bound for $\sup_{x\in \mathbb{R}_+^0}\abs{f(x)^{\gamma-1}}$ in terms of the norm of $f$, and we can not bound $\norm{f^\gamma}$ in terms of $\norm{f}$. This prevents us from obtaining locally bounded and linear growth properties for $0<\gamma < 1$ even when restricting to functions in $\h_>$. 
\end{remark}

We now want to show the locally Lipschitz property for 
$f\mapsto f^\gamma$ and $f\mapsto\mathcal{M}_{f^{\gamma}}$ by applying, respectively, Corollary \ref{lem:specialcase:mult:structure} and Corollary \ref{lem:specialcase:mult:structure2}. Then we must show that $\phi^\gamma(y)=y^{\gamma}$ has locally Lipschitz derivative $\left(\phi^{\gamma}\right)^\prime(y)=\gamma y^{\gamma-1}$. However, the function $\left(\phi^{\gamma}\right)^\prime(y)$ is itself differentiable with derivative $\left(\phi^{\gamma}\right)^{\prime\prime}(y) =\gamma(\gamma-1) y^{\gamma-2}$ which is unbounded  for $1<\gamma<2$ on any interval of the form $(0,\delta)$, for $\delta >0$. Therefore, we can apply Corollary \ref{lem:specialcase:mult:structure} and Corollary \ref{lem:specialcase:mult:structure2} only for $\gamma\in\N$ and for $\gamma \in \mathbb{R}$ with $\gamma \geq 2$. We shall deal with the case  $1< \gamma < 2$ in Proposition \ref{prop:gamma12} by introducing a state-dependent exponent $\gamma$ for the CEV-SPDE.

\begin{cor}\label{cor:lipgamma}
		The following statements hold:
		\begin{enumerate}
			\item Let $\gamma \in \mathbb{N}$. Then $f\mapsto f^\gamma$ and $f\mapsto \mathcal{M}_{f^{\gamma}}$ are locally Lipschitz  maps in $\h$.
			\item Let $\gamma \in \mathbb{R}$ with $\gamma \geq 2$.  Then $f\mapsto f^\gamma$ and $f\mapsto \mathcal{M}_{f^{\gamma}}$ are locally Lipschitz  maps both in $\h_+^0$ and in $\h_+$.
		\end{enumerate} 
	\end{cor}
	\begin{proof}
		By means of Corollary \ref{lem:specialcase:mult:structure}, we must show in the two cases that the function $\phi^\gamma(y)=y^{\gamma}$ 
		has locally Lipschitz derivative. For $\gamma=1$, the derivative $\left(\phi^{\gamma}\right)^\prime(y)=1$ corresponds to the constant map, which is clearly Lipschitz. Let now $\gamma\in\N$ with $\gamma \geq 2$. By the mean value theorem, we know that for every $y_1, y_2\in\R$ we can write $\left(\phi^{\gamma}\right)^\prime(y_1) - \left(\phi^{\gamma}\right)^\prime(y_2)= \left(\phi^{\gamma}\right)^{\prime\prime}(c) (y_1-y_2)$ for some $c\in [y_1,y_2]$. This implies that for every $n>0$, $\abs{\left(\phi^{\gamma}\right)^\prime(y_1) - \left(\phi^{\gamma}\right)^\prime(y_2)}\leq \gamma(\gamma-1) n^{\gamma-2} (y_1-y_2)$ for every  $y_1,y_2 \leq n$, hence $\phi^\gamma$ has locally Lipschitz derivative. The same holds for $\gamma \in \mathbb{R}$ with $\gamma \geq 2$ as long as $y\geq 0$. This shows, by means of  Corollary \ref{lem:specialcase:mult:structure}, that $f\mapsto f^\gamma $ is locally Lipschitz. We further conclude from  Corollary \ref{lem:specialcase:mult:structure2} that also $f\mapsto \mathcal{M}_{f^{\gamma}}$ is locally Lipschitz in both cases.
	\end{proof}

We point out that Corollary \ref{lem:specialcase:mult:structure} and Corollary \ref{lem:specialcase:mult:structure2} provide conditions that ensure both locally bounded and locally Lipschitz properties. However, we have applied them in Corollary \ref{cor:lipgamma} only to show the local Lipschitz property. Instead, to show the local boundedness in Corollary \ref{lem:fgamma2} we applied Proposition \ref{prop:locallybounded} which refers to a general $\Psi$, namely not exploiting that $\Psi^\gamma(t, f)$ is of the form of \eqref{eq:betaphi}. Indeed, as we can see from Corollary \ref{lem:fgamma2} and Corollary \ref{cor:lipgamma}, we have the local Lipschitz property  only for $\gamma\ge 2$, while we have locally bounded also for $1\le \gamma <2$. Hence, since locally Lipschitz and locally bounded properties require different restrictions on the exponent $\gamma$, we have to deal with them separately.

It might be surprising at first that the map $f\mapsto f^\gamma $ is not locally Lipschitz on $\h_+^0$ for $1<\gamma <2$, even if the function $y\mapsto y^\gamma$ is locally Lipschitz on $\R_+^0$. This, as for Corollary \ref{lem:fgamma2},  is due to the nature of the Filipovi\'c norm, from which we need to restrict ourselves to $\gamma\ge2$. Indeed, by replacing $f^\gamma$  with $f^{\gamma-1}$ in equation \eqref{normfgamma},  we see that for $\gamma < 2$ the term $f(x)^{\gamma-2}$ gets large for those values of $x$ for which $f(x)$ is small. Hence $\norm{f^{\gamma-1}}_{}$ might not be bounded, and we can not guarantee the local Lipschitz continuity of $f\mapsto f^\gamma$ in the Filipovi\'c space. Notice that, similarly to Remark \ref{remark:noboundf}, we can not guarantee locally Lipschitz even if restricting the domain to $\h_>$.
	
	However, it is possible to overcome this issue by changing the behavior of the function $\phi^\gamma$ in $0$. This can be achieved by making the exponent $\gamma$ state dependent.

	\begin{prop}\label{prop:gamma12}
		Let $\gamma >1 $ and let $\tilde{\gamma}: \mathbb{R}_+^0 \rightarrow [1,\infty ) $ with bounded derivative on $\R_+$ be such that:
		\begin{enumerate}
			\item $\lim_{y\rightarrow \infty } \tilde{\gamma } (y)=\gamma$;
			\item There exists $\varepsilon >0$ such that $\tilde{\gamma } (y)=1$ for all  $y\in [0,\varepsilon]$.
		\end{enumerate} 
		Then for $\mathcal{Y}\in\{\h_+^0, \h_+\}$  the map
		\begin{align*}
			\Phi^{\tilde\gamma}: \mathcal{Y} &\to\mathcal{Y} \\
			f &\mapsto \{ x\mapsto \phi^{\tilde\gamma}(f(x)):=f(x)^{\tilde{\gamma}(f(x))}\}
		\end{align*}
		is well defined and locally Lipschitz. 
	\end{prop}
	\begin{proof} We first notice that the function $\tilde\gamma$ is locally bounded because of assumptions \emph{1.} and \emph{2.},  and locally Lipschitz because of the assumption of bounded derivative. Then, for every $n$, there exists $M_n>0$ such that $|\tilde\gamma(y)|\le M_n$ for all $y\leq K_\delta n$, 
	and 
	\begin{equation*}
		|\phi^{\tilde\gamma}(y)| = |y^{\tilde\gamma(y)}| \le n^{M_n} <\infty,
	\end{equation*}
	which shows that $\phi^{\tilde\gamma}$ is locally bounded. Similarly, it follows that $y^{\tg (y) -1}\leq \max \{1 ,n^{M_n - 1} \} $. Moreover, we notice that $\tg '(y )=0$ for $y\in [0,\varepsilon )$ because of the assumption $\tg (y )=1$ for $y\in [0,\varepsilon]$. Hence  we deduce that there exists an $\bar{M}_n>0$ such that
	\begin{equation}
		\label{eq:derphigamma}
		\left(\phi^{\tilde\gamma}\right)^\prime(y) =  \tg (y)  y^{\tg (y) -1} + y^{\tg (y)} \log (y) \tg ' (y )\le \bar{M}_n,
	\end{equation}
	for $y\leq K_\delta n$. The map $\phi^{\tilde\gamma}$ is therefore locally Lipschitz with Lipschitz constant $\bar{M}_n$ for $y\leq K_\delta n$. 
	
	From the equality in \eqref{eq:derphigamma} we also obtain that $\left(\phi^{\tilde\gamma}\right)^\prime(f(x))$ is bounded, since $\lim_{x\rightarrow \infty } f(x)< \infty$. Moreover, with $\left(\Phi^{\tilde\gamma}\right)^\prime(f)(x) = f' (x) \left(\phi^{\tilde\gamma}\right)^\prime(f(x))$, we obtain that
		\begin{equation*}
			\| \Phi^{\tilde\gamma}(f)\|^2 = \left(\phi^{\tilde\gamma} (f(0))\right)^2 + \int_0^{\infty} \left(f' (x) \left(\phi^{\tilde\gamma}\right)^\prime(f (x)) \right)^2 w (x) dx <\infty,
		\end{equation*}
		which shows that $\Phi^{\tilde\gamma}:\mathcal{Y}\rightarrow\mathcal{Y}$ is well defined.
		
		We now prove that $\Phi^{\tilde\gamma}$ is locally Lipschitz. Let $f_1, f_2 \in \mathcal{Y}$ and $\tilde M_n:=\max \{M_n,\bar{M}_n\}$. Then
		\begin{equation*}	
			\begin{aligned}
			\| \Phi^{\tilde\gamma} (f_1) - \Phi^{\tilde\gamma}(f_2) \|^2 =& \left( \phi^{\tilde\gamma}(f_1(0)) - \phi^{\tilde\gamma}(f_2(0)) \right)^2  \\&+ \int_0^{\infty} \left(f_1' (x) \left(\phi^{\tilde\gamma}\right)^\prime( f_1 (x)) - f_2' (x) \left(\phi^{\tilde\gamma}\right)^\prime( f_2 (x)) \right)^2 w (x) dx.
			\end{aligned}
		\end{equation*}
	In particular, since $\phi^{\tilde\gamma}$ has Lipschitz constant bounded by $\tilde M_n$ on $[0,K_\delta n]$, for the first term on the right-hand side of the equation above, we observe that $\left( \phi ^{\tilde\gamma}(f_1(0)) - \phi^{\tilde\gamma}(f_2(0)) \right)^2\leq K_\delta^2\tilde M_n^2  \| f_1 - f_2 \|^2 $ as long as $ \| f_1 \|,  \| f_2 \|\leq n$. Moreover, for  the second term we can write that
  \begin{equation}
      \label{eq:proofsigma2}
		\begin{aligned}
			&\int_0^{\infty} \left(f_1' (x) \left(\phi^{\tilde\gamma}\right)^\prime( f_1 (x)) - f_2' (x) \left(\phi^{\tilde\gamma}\right)^\prime( f_2 (x)) \right)^2 w (x) dx	 \\
			&\leq  2 \int_0^{\infty} \left( f_1' (x)  - f_2' (x) \right)^2 \left(\phi^{\tilde\gamma}\right)^\prime( f_1 (x) )^2 w (x) dx \\&\quad\,+ 2 \int_0^{\infty} f_2' (x)^2 \left(\left(\phi^{\tilde\gamma}\right)^\prime( f_1 (x)) - \left(\phi^{\tilde\gamma}\right)^\prime( f_2 (x)) \right)^2 w (x) dx .
		\end{aligned}
   \end{equation}
		In particular, $f_1(x)\leq K_\delta \| f_1 \|$, hence here $\left(\phi^{\tilde\gamma}\right)^\prime(f_1(x)) \le  \tilde M_n $ for any $ \| f_1 \| \leq n$.  The first term in \eqref{eq:proofsigma2} is therefore bounded by $2 \tilde M_n^2 \| f_1 - f_2 \|^2 $,  for any $f_1,f_2$ with $\| f_1 \|,\| f_2 \|\leq n$. Moreover, by Lipschitz continuity of  $\left(\phi^{\tilde\gamma}\right)^\prime$ on $[0,K_\delta n]$, the second term on the right-hand side of \eqref{eq:proofsigma2} is bounded by $2 n^2\tilde M_n^2 K_\delta^2\| f_1 - f_2 \|^2$, as long as $\| f_1 \|, \| f_2 \| \leq n$. Combining all together, we get that
		\begin{equation*}
			\| \Phi^{\tilde\gamma} (f_1) - \Phi^{\tilde\gamma}(f_2) \|^2 \le 2 \tilde M_n^2\left( 1+n^2K_\delta^2\right)\norm{f_1-f_2}^2,
		\end{equation*}
		which concludes the proof. 
	\end{proof}

We conclude this study on the CEV-SPDE with the following corollary, showing that, by means of Corollary \ref{cor:positivity_betaphi} and Lemma \ref{lem:positivity}, we can state conditions ensuring the positivity of the solutions of the CEV-SPDE.
	\begin{cor}
		Let $\beta:\R_{+}^0\to \R_+$ be differentiable such that  $\inf_{t\in\R_{+}^0}\frac{\beta'(t)}{\beta(t)}>-\infty$, and let $\phi^\gamma:\R_+^0\to\R_+^0$ be defined by $\phi^\gamma(y):=y^\gamma$. If $\gamma \ge 2$, then for every $g_0\in\h_+$, respectively for every $g_0\in\h_>$, the local mild solution $(g, \tau)$ to the CEV-SPDE \eqref{spde_cev}  with zero drift belongs to $\h_+$, respectively to $\h_>$ for every $t\le \tau$.

        Let further $\alpha:\R_+^0\times \h_+\to \h$ satisfy locally bounded and locally Lipschitz conditions accordingly to the results in Section \ref{sec:conditions}, and for any initial condition $\bar g_0\in \h_>$ denote by $(\bar g , \bar\tau)$ the local mild solution of the SPDE \eqref{spde_cev} with drift $\alpha$, where $\bar\tau\le \tau$ almost surely. If condition \eqref{novikov:cond_positivity}	is satisfied for some $\bar{T}<\tau$, then $\bar g_{t}\in \h_>$ for any $t<\bar T$.
	\end{cor}
	\begin{proof}
		It is easy to see that condition \emph{1.} of Corollary \ref{cor:positivity_betaphi} is satisfied. For condition \emph{3.} we need to show that 
		$\inf_{y\in (0,\varepsilon)}y^{\gamma-1}>-\infty$, which holds for $\gamma\ge1$. However, in order to guarantee existence and uniqueness of the solutions of the CEV-SPDE, we must restrict ourselves to $\gamma \ge 2$ accordingly to Corollary \ref{cor:lipgamma}. This gives us positivity of the solution for the CEV-SPDE with zero drift. The second part of the statement is a consequence of Lemma \ref{lem:positivity}.
	\end{proof}
In this section, we have considered examples that highlight the links between the Hilbert space-valued approach with point-wise operating coefficients and one-dimensional modeling of the forward for a fixed delivery time $T$. In our forthcoming paper \cite{DeteringLavagnini}, we propose a specification that resembles the local volatility model used in equity markets \cite{dupire}, and we perform a calibration to market data.

\section{Conclusions and outlook}
\label{sec:conclusions}
We introduced and studied a class of Hilbert-space valued state-dependent models for HJM forward curves, defined via point-wise operating maps and multiplicative operators. The advantage of our approach is twofold. On the one hand, thanks to the SPDE and Hilbert space approach, the entire forward curve is modeled with one dynamics. This allows for a coherent state space, for a control of the shape of the forward curves, and for the flexible modeling of all possible forwards without arbitrage contradictions across different delivery times. On the other hand, with our proposed specification based on a point-wise operating structure of the coefficient functions, the mathematical and numerical complexity which is usually inherited from the Hilbert-space approach can be reduced. This is thanks to a clear tie to classical one-dimensional forward models. In fact, our specification significantly reduces numerical complexity if a derivative on a fixed maturity contract is to be priced. In particular, one can either resort to analytical formulas or simulation and discretization schemes for one-dimensional SDEs. This makes the proposed approach flexible enough to define state-dependent models that can be calibrated to the market. Furthermore,  it has the potential to make high-dimensional models more attractive to practitioners, who still need to benchmark to simple models, such as, e.g., the Black-Scholes model. In the accompanying paper \cite{DeteringLavagnini}, we shall use our class of state-dependent models in order to construct a local-volatility framework, in the spirit of \cite{dupire}, for commodity and energy markets, with the aim of fitting the volatility smile that is typical for these markets.

Our work gives rise to interesting future research. The maps that we have considered here induce a class of locally state-dependent models. A possible extension is to consider instead maps $\Psi:\h\to\h$ of the form $\Psi(f)(x) = \psi(f(x), f'(x))$, for some $\psi:\R\times \R\to\R$, with $f\in\h$ and $x\in\R_+^0$. That is, one could consider maps that act on both the current level of the curve and its slope. This particular specification may potentially lead to more accurate models as the state dependency would be linked to the level of forwards with similar maturities. However, it requires one to work with an extension of the Filipovi\'c space to twice (weakly) differentiable functions, such as the spaces proposed by \cite{FILWIL18}, \cite{filipovic22a} and \cite{filipovic22b}. Furthermore, maps may be considered that act directly on the time to maturity, namely $\Psi:\h\to\h$ of the form $\Psi(f)(x) = \psi(x, f(x))$, for some $\psi:\R_+^0\times \R\to\R$, with $f\in\h$ and $x\in\R_+^0$. This specification would be in line with Remark \ref{remark:psix} and would allow, for example, to include a maturity effect, such as the Samuelson effect \cite{samuelson}, directly in the SPDE coefficients. Note that in the current setting, such an effect can already be incorporated through the noise structure.

\subsection*{Acknowledgment}
The authors are grateful to Fred Espen Benth, Paul Kr\"{u}hner, and Sara Svaluto-Ferro for fruitful discussions.

\bibliographystyle{plain}  
\bibliography{bib}
\end{document}